\documentclass[12pt]
{amsart}
\usepackage{amsfonts,amsmath,amssymb,amsthm,graphicx,afterpage,url}
\usepackage{hyperref}
\hypersetup{
   colorlinks   = true, %Colours links instead of ugly boxes
   urlcolor     = blue, %Colour for external hyperlinks
   linkcolor    = blue, %Colour of internal links
   citecolor   = red %Colour of citations
}
\usepackage{color}

\input{xy}
\xyoption{all}

\def\R{{\mathbb R}} \def\Z{{\mathbb Z}}

\long\def\comment#1\endcomment{}

\newcommand{\aronly}[1]{}

\theoremstyle{theorem}
  \newtheorem{theorem}{Theorem}[section]

\theoremstyle{definition}
  \newtheorem{remark}[theorem]{Remark}
  \newtheorem{example}[theorem]{Example}

\begin{document}

\title{On different reliability standards \\ in current mathematical research}

\author{A. Skopenkov}

%https://retractionwatch.com/

%use shortened versions from bit.ly or tinyurl.com
%combined with \usepackage{url} or \usepackage{hyperref} to make them clickable

%Joji\' c, D. et al. The coloured Tverberg theorem, extensions and new results Izvestiya: Mathematics. 86:2 (2022)
%https://arxiv.org/pdf/2108.08804 \S4

\thanks{Moscow Institute of Physics and Technology, and Independent University of Moscow.
Email: \texttt{skopenko@mccme.ru}.
\texttt{https://users.mccme.ru/skopenko/}.
\newline
I am grateful to R. Karasev and M. Skopenkov for useful discussions, and to A. Sossinsky for style editing of a part of the text.
\newline
Key words and phraces: research integrity, peer review journals, scientific ethics, history of mathematics.
}
%\newline
%This version of arXiv:2101.03745 is updated more often than the arXiv version.}

\date{}

\abstract In this note I describe reliability standards for writing and reviewing mathematical papers; these standards are (in my opinion) vital for the progress of mathematics.
I give examples of applying the described or other reliability standards.
%This is a preliminary version, please do not distribute.
\endabstract

\maketitle

\tableofcontents

\hfill{\it Mozart was a great composer.}

\hfill{The Barber of Siberia}
%,  N. Mikhalkov

\section{Introduction}\label{s:intr}

The reader knows the difference between a medical drug that passed the tests and is acceptable for use,
and a drug that did not.
It is vital to have the tests, the information on passing them, and the information on how competent the tests were.
The purpose of this paper is to attract attention to the analogous problem of reliability in mathematics.

This paper will hopefully be interesting not only to mathematicians and scientists, but also to people interested in mathematics as a part of culture, and as an activity in which their taxes are spent.
Most of this paper is  accessible to non-specialists (although some examples do contain technical details).

By `reliability' of a research paper I understand the presentation of its results and arguments ensuring their validity.\footnote{\label{f:clar} Thus reliability is not the same as clarity, but is closely related and correlated to clarity  (of the results, of their proofs, and of relations to earlier known results).
For recommendations on clear writing see the references in \cite[\S2.1]{Pa18}, and \cite{Pa18}.
Other obvious requirement from a research paper is novelty of results.}
%This is independent of the value and importance of results.
Reliability standards in mathematics change.
E.g. a long time ago, oral proof was considered reliable.
As proofs gradually became more complicated, there appeared the understanding (and then the rule) that only written proof can be considered to be reliable or not (and so could constitute a fair claim for a result).

\smallskip
{\it `All times are changing times, but ours is one of massive, rapid moral and mental transformation.
Archetypes turn into millstones, large simplicities get complicated, chaos
becomes elegant, and what everybody knows is true turns out to be what some people used to think.
It's unsettling.
For all our delight in the impermanent, the entrancing flicker of electronics,
we also long for the unalterable... Don Quixote sets out forever to kill a windmill.}

(U. K. Le Guin, Tales From Earthsea)

\smallskip
Different reliability standards coexist.
So it would be nice if mathematicians and math journals publicly and explicitly revealed their reliability standards  by publishing their (potentially different) opinions on specific examples.
This would allow more competent decisions of the math community, sponsors, and tax-payers to support this trend or another in mathematics.

%(The same applies to refereeing standards, see below.)

%`Our diversity is our strength.' (A. Jolie)

\smallskip
{\it It is not merely true that a creed unites men.
Nay, a difference of creed unites men --- so long as it is a clear difference.

I am quite ready to respect another man's faith; but it is too much to ask
that I should respect his doubt, his worldly hesitations and fictions, his
political bargain and make-believe.}

(G. K. Chesterton, What's Wrong With The World)

%It would be helpful if mathematicians will publish short expositions of their reliability standards
%%(or to send their complements to this paper, to be published without change),
%(no matter whether their reliability standards are different from mine or not).

\smallskip
A {\bf reliable reference} is a text whose results mathematicians can use without doing the author's work:
to check thoroughly the validity of the results and to properly describe the relationship to earlier publications.\footnote{We should at least try to understand arguments even from unreliable references.   
However, to understand an argument omitting details on the assumption that they are already checked requires much less effort than to check an argument including all details on the understanding that every detail might potentially be a problem.
Everybody who uses a computer program
%(or a water tap)
knows the big difference between a computer program and a computer program that works.
Everybody who performs work for a user knows that quite some time
is required from the developer to test the result.}
Being a reliable reference in the first approximation corresponds to `accept as is' recommendation for a peer reviewed journal.\footnote{Such a recommendation allows further changes that need not be checked by a referee (unless the changes are significant).}
The peer review system \cite{Pe} is reviewed in \S\ref{s:peer}--\S\ref{s:isrask17}, see also \cite{Paj} and the references therein.

%??? Formal discussions not informal

The reliability standards I share are presented in \S\ref{s:differ} and \S\ref{s:peer} below mostly in the form of recommendations. These recommendations are based on traditions I learned, as well as
on my own experience as an author, a reader, a referee, and a jury member of scientific contests.

%\footnote{Some readers shared their feeling that the purpose of these examples is to defend my
%priority for some results, in spit of the opposite being clear from the text, see also footnote \ref{f:help}.}

This paper can be considered as a complement of practical nature to the discussion by prominent mathematicians \cite{JQ93, ABC+}.
For similar issues on references, on computer science conferences, and on the $abc$-conjecture see
\cite{Pan, Fo, Kl18}, respectively.

This paper can make only a modest contribution to revealing and improving of standards for writing,
discussing, refereeing, and taking editorial decisions.
However, collective work of thousands of researchers on their usual tasks (like writing a paper, a referee report or a Zentralblatt/MR review) can make mathematical research reliable according to the moral and material
support that we are applying for.

\begin{remark}\label{r:exam}
General words on delicate subjects (like how to play piano) are almost useless without examples.
It is only by going into substantial details that we can conclude if a referee report (or an author's reply to such a report) is competent; otherwise we are likely to be distraught by emotions.

Also, it is easy to violate our declared principles, when they come into contradiction with our worldly
hesitations and fictions, group thinking \cite{Gr}, political bargain and make-believe.
(This is usually done unintentionally and unconsciously, and revelation of this is painful to us.)
The turning point is to practice what one preaches, and this inevitably changes the preaching.
{\it `He who practises the Way (T$\hat a$o) suffers daily loss of its false shine'.}\footnote{This is an English translation of a citation from the Russian translation \cite[Chapter 22]{Ch} of {\it Chuang tzu}.
I am grateful to Yan Pan for informing me that the following translation (said to be by Y. Lin) seems more acceptable  to Chinese scholars: `He who practises the T$\hat a$o, daily diminishes his doing'.
In particular, the words `of its false shine' are not present either in {\it Chuang tzu} or in English or Russian translations different from \cite{Ch} and available to me.
Still, because of the preceding text in \cite[Chapter 22]{Ch} I find `of its false shine' a proper commentary.}

%and I hope Chuang tzu had something similar in mind, although he had not this in his writing.
%are not present in \cite[Chapter 22]{Ch} but are present in Chinese-Russian translation.
%It would be nice to have an opinion of a specialist in Chinese culture on whether the above
%longer translation is more authentic, cf. \cite[Chapter 22]{Ch}.

Thus examples are important.
They are presented in \S\ref{s:differex}, \S\ref{s:peerex}, \S\ref{s:isrask17}, \cite{Ha},
\cite[Remark 5.3, \S5.3, \S5.4]{Sk16}; references to examples are presented in \S\ref{s:differ} and \S\ref{s:peer};  see also \url{https://scirev.org/}.
\end{remark}

{\bf Conventions.}
I omit `in my opinion' for brevity.
Some {\it general} statements (`should be', `clearly', etc.) are made
%`I consider', `here is a procedure',  means that I do not expect
on the assumptions that no public objections will appear; I am ready to give justification, if public objections appear.
%As opposed to that,
All {\it specific} examples are justified as thoroughly as to possibly annoy some readers with details; such readers can omit the details.

{\it A genius makes his own rules, but a `how to' article is written by one ordinary mortal for the benefit of another... Authors of articles such as this one know that, but in the first approximation they must ignore it, or nothing would ever get done.} \cite{Ha74}

By {\it user} I mean the user (reader, reviewer, etc.) of a paper, who can also be the developer
(author, advisor, etc.) of other papers.
Often users read papers for
%the purpose of
developing other papers.

%\newpage
\section{Different reliability standards}\label{s:differ}

\begin{remark}\label{r:frompt}
(a) Among the first steps of checking the proof it is advisable to

\quad (1) write on the first pages the statements of the results, together with all the definitions used;

\quad (2) structure the paper so as to separate motivations from statements and proofs,
e.g. by having a separate subsection `motivations' (see e.g. \cite{Sk20o});

\quad (3) give rigorous proofs of all statements;\footnote{This can be done either under a head `proof of such and such numbered statement', or in phrases like `Such and such numbered statement follows from [a list of  numbered statements from this or other papers]', `The proof is obtained from such and such a proof [reference] by the following substitution: [replace this by that, etc.]'.
\newline
This should be done instead of, or in addition to, informal phrases like `it is somewhat well-known to be', or instead of some informal observations ended with `we obtain the following: [statement]' (see e.g. \cite[the third paragraph of footnote 4]{Sk20e}).
Phrases like `The proof is analogous to such and such a proof [reference]' should be avoided when the analogy requires too much time and efforts from a user.}

\quad (4) study proofs of most closely related results and send the paper to their authors for comments (see more in Remark \ref{r:steps}.a).

Using unsophisticated language, and sticking to statements and proofs of the main results,
helps to find (and correct) mistakes.

%looking from a different point

(b) It is advisable to explicitly state any non-trivial result, and either prove it or give a reference.
This is not only necessary for reliability, but contributes to the unity of mathematics,
by making a paper more accessible to mathematicians from other areas.
This in turn ensures higher reliability.
% (in particular, in checking novelty).
Lack of explicit statements decreases reliability, and also contributes to artificial splitting of mathematics (and even of its areas) into different subjects whose representatives cannot use each other's work.

(c) The usual problem related to the lack of explicitly given rigorous definitions / statements is not that it is hard to reconstruct them, but that this can be done in several ways.
\end{remark}

\begin{remark}\label{r:gener}
(a) {\it Listeners} [AS: and readers] {\it are prepared to accept unstated (but hinted) generalizations much more than they are able ... to decode a precisely stated abstraction and to re-invent the special cases that motivated it in the first place.} \cite{Ha74}

(b) {\it The modern world is full of theories which are proliferating at a wrong level of generality, we're so good}  at theorizing, {\it and one theory spawns another, there's a whole industry of abstract activity which people mistake for thinking.} (I. Murdoch, The Good Apprentice)

(c) See \cite[Preface]{GKP} for discussion of similar issues.

(d) It is advisable to postpone technical results and definitions to later (sub)sections and bring bright results to earlier (sub)sections.
In my opinion, bright results and their proofs are potentially more useful than technical versions of known constructions which so far did not yield any bright results.
Here by bright results I mean non-trivial results whose statements are accessible to mathematicians specialized in this area of mathematics, but not necessarily in the subject of the paper.

(e) When one uses a specific theory, it is advisable to explicitly state results to be proved with the help of this theory, but in terms not involving the theory (see e.g. \cite{Sk16, Sk19, Sk20e}).
This makes the {\it application} of the result  accessible to mathematicians who have not specialized in the theory.
So one is motivated to study the theory and sees explicit statements which could guide this study.

%, so an exposition shaved by Occam's razor
%(and so accessible to non-specialists) is hopefully interesting.

Instead of the above way, some papers start exposition with details of a specific theory which are matter-of-fact to specialists but are inaccessible to mathematicians from other areas.
(E.g. compare \cite{LL18} to \cite{Sk19} and \cite[\S3.3]{BZ16} to \cite[\S2.3]{Sk16}; see \cite[Remark 5.7.8]{Sk16}.)
Sometimes this happens because of the false (and possibly not conscious) assumption that the readers will accept artificial sophistication as depth and high non-triviality.

(f) Sometimes originally a complicated proof is invented.
%(e.g. 30-page solution of a quadratic equation using Galois theory).
(This often happens because of lack of familiarity with simple expositions of the subject.)
Freeing a proof from complications appearing in its invention
%(but not necessary for the proof)
is a way of checking the proof, not only the courtesy of presenting a simplified proof.
%However, a reader wants a short proof not made artificially inaccessible by making it dependent on more %complicated ideas unnecessary for the short proof (even if those ideas did appear in the authors' approach).

In particular, a large part of a theory that authors describe could be superfluous for a proof of the result.
%A contradiction is at the end an equality between different numbers
%(or another simple and basic objects like numbers).
It is advisable to give a direct proof without formally referring to all the theory that helped the authors to invent the proof.
Relation to the theory can though be mentioned in a remark which is not part of the proof.
See also Remark \ref{r:gener}.e.

Analogous remarks holds with `proof' replaced by `statement', etc.

(E.g. compare longer proofs in \cite[\S3,\S4]{KS20}, \cite{Pa15}, \cite{Pa20} to shorter proof in \cite[\S2]{KS20}, \cite{Sk18o}, \cite{PS20}, respectively, or proofs in different versions of \cite{Pan15, AMSW}.)
\end{remark}

\begin{remark}\label{r:fromll}\footnote{For example of a different opinion by L. Vok\v r\'inek see \cite[footnote 4]{Sk20e}}
(a) A published paper is for a much wider audience of mathematicians than just referees, close colleagues and Editors.
So motivations for main results, details of the proof, etc. should be written in the paper, not in letters to referees, close colleagues and Editors.

(b) A published paper is for users, not for developers.
Working on details could be an interesting task for a developer, but is usually not within the intents of a user.

The following are good lower estimates of how hard work on details is:

$\bullet$ the amount of time required for authors
(or for other mathematicians) to make the details publicly available upon request of a reviewer.

$\bullet$ the amount of text written by the authors to justify that the details need not be made publicly available.

%text (i.e. in arXiv or in a journal),
%and to reveal real or imaginary flaws in the correspondence asking for details.

(c) Updating the arXiv version is not considered as an indication that the published version has any serious gaps.
Some authors previously updated arXiv versions upon my suggestions as a math reviewer, and we didn't have any discussion about that.
See e.g. arXiv:1609.06573v3,   arXiv:1209.1170v4 and a forthcoming paper by D. Gugnin.
%A less positive example is arXiv:1512.05164v6 +  arXiv:1808.08363v2.

(d) Discussions of a text involving a user of this text should refer (at least upon request of the user) to the text, not to any non-existent text obtained from the discussed text by some changes (see e.g. \cite[the third paragraph of footnote 4]{Sk20e}).
Of course, this need not apply to discussions of a text between developers (e.g. coauthors) of this text.
\end{remark}

\begin{remark}[important steps to prepare a quality submission]\label{r:steps}
(a) The mathematical community needs the reliability of a paper to be checked {\it before} not {\it after} putting the paper to arXiv (see motivations in Remark \ref{r:clare}).
The same is true for {\it preliminary} checking of the novelty, because putting a paper to arXiv is a means of {\it further} checking of the novelty.
So it is advisable before putting a paper to arXiv to discuss it among specialists in its area.
Such a pre-submission discussion usually allows the author to check whether his/her
results are clearly stated, new, and completely proved.
This allows to improve quality of the paper.

Improving quality may involve a disillusion and dissatisfaction (as a part of learning, cf. Remark \ref{r:exam}).
A dissatisfaction which might appear during such work is a natural part of improving quality of the paper and qualifications of the author.
Such work is interesting if authors (=developers) recognize the importance of learning and fulfilling the  wishes of their readers (=users).
Such work is annoying only if authors write under the assumption (however unconscious) that their work need not be useful.
Improved quality of the paper publicly available on arXiv improves the author's reputation, while low quality damages it.

It is advisable to put a paper on arXiv before submitting it to a journal, cf. Example \ref{r:review}.
Without this simple procedure the possibility remains that the results of a paper are already (partly) known.
What is unknown to one group of mathematicians can be known to another group.
Putting a paper on arXiv allows including into a pre-submission discussion (described in the previous paragraphs) people who work in related  areas, but are not in contact with the author.
Improved quality of a paper published in a journal improves the author's reputation, while low quality damages it
even more significantly.

In a less formal pre-submission discussion it is easier to help the author, to share ideas with him/her,
and to minimize the critical part of such help.
Publication of a paper on arXiv without prior discussion with a colleague means that the author expects a public,  not a private approval or criticism of this colleague.
The colleague might still prefer private criticism, but could be compelled to criticize publicly if the arXiv
text contain flaws which obstruct the progress of mathematics, cf. Remark \ref{r:clare}.
%are not corrected for a long time.

Writing a referee report on a paper is a responsible task involving double-checking (see Remark \ref{r:qual}).
In this time-consuming form it is much harder to help the author than via informal discussions.
%see e.g. Remark \ref{r:histjo}.

The above steps do not absolutely protect against significant flaws, see e.g. \cite[footnote 3 in p. 2]{Sk08}.
However, with the above steps done, the responsibility is shared with the math community.

(b) During pre-submission discussion, specialists in the area might send their specific suggestions/criticism which they consider important (below this is shortened to just `suggestions').
Then it is advisable to put on arXiv (or submit to a journal) a revised version approved by specialists.
Of course this might not be easy to do.
E.g. the authors can receive no comments from somebody; they can
disagree with some suggestions (and they cannot be sure that they would not receive another stupid or essential suggestions the day after they finalized their work, based on previous suggestions).
Hence the authors can decide to submit their paper to arXiv/journal even if they did not

\quad (1) receive a feedback from some persons,

\quad (2) take into account some suggestions, or

\quad (3) give a chance to a person who sent feedback to learn his/her opinion whether his/her
critical remarks are properly taken into account.
%suggestions are properly incorporated.

There is no need to mention (1), but it is advisable to mention (2) or (3) in the submitted text.
See Example \ref{r:stepsex}.
\end{remark}

%So this work should be done
%by developers (coauthors, advisors, etc) , and should involve discussions between developers.
%before arXiv publication, journal  publication and defense of Ph. D. thesis.
%(or, in the worst case, in the first arXiv version following the authors' work on a counterexample
%to published version).

\begin{remark}[claims and responsibility]\label{r:clare}
Here I present some motivations for the recommendations of Remark \ref{r:steps}.
A natural reaction to Remark \ref{r:steps} is as follows:
{\it The opinions of people vary on how to use arXiv.
It is good if there is some diversity.
Some people would put on arXiv only a finished paper which appeared in a journal,
some people would put there very early preprint.
Overall, a paper on arXiv may be both: already a solid work or something just very early.
Many mathematicians, if they find some problems in an early preprint on arXiv, would either ignore this, or
send their remarks to the authors.}
%(There are also very subtle differences between the options (i) there is a major mistake in the paper
%(ii) the presentation is very sloppy (iii) I do not like the presentation but perhaps some people
%with enough effort can understand it.)

A research paper (on arXiv or elsewhere) can have both positive and {\it negative} impact
on the development of mathematics.
A usual example of a negative impact is as follows.
A paper claims an interesting result but the proof does not form
%is not complete enough to provide
a reliable reference (see \S\ref{s:intr}).
Then this result cannot be used.
Also, colleagues are discouraged from providing a reliable reference because its publication might involve hard and unpleasant work of justifying that the existing text is not a reliable reference.
So the paper obstructs the appearance of a reliable reference (by authors of the paper or by other authors).
So the paper may involve unfair competition (cf. footnote \ref{f:help}).

Mathematicians often have interesting ideas but no time to develop them into reliable references.
Sharing preliminary ideas in the form of claiming results has a negative impact described above, but could also have positive impact.
Sharing preliminary ideas as preliminary ideas is useful, not harmful.

A \emph{clear} diversity in arXiv papers, when a paper in its first lines explains  whether the authors think it is a finished paper or a very early preprint, certainly stimulates progress of mathematics.
%See e.g. \cite{Sk19},
%\cite{Sm} (personally, I would prefer even a clearer explanation of what is final and what is not in that
%long paper containing many interesting ideas),
%more examples???.
Lack of clear diversity allows the attempt, however unintentional, to both make a claim and not have the responsibility for making it.
With some administrative support, this attempt can well be successful in terms of having paper published, getting a grant or a job, etc.
So without \emph{clarity} (provided either by authors or by their critics) the style of not caring for users has a significant advantage of spreading comparative to the style of caring for users.
This has negative impact on the development of mathematics.
\end{remark}

%Due to efforts of more careful authors arXiv papers are usually believed (at least if written
%This motivates recommendations of Remark \ref{r:steps}.

\begin{remark}\label{r:prework}
This remarks concerns (not reliability but) clarity of exposing relations to previously known results, see footnote \ref{f:clar}.

It is necessary to cite previous research, even if it is not a formal part of the presented results.
% from the current paper.
E.g. if previous research is closely related (see e.g. \cite[footnote 3 in p. 2]{Sk08}), or introduces an important idea of the current paper (see Example \ref{e:pt19}), or presents an alternative exposition which has a chance to be simpler for some readers (see Example \ref{e:pt19}).
\end{remark}

\section{Different reliability standards: examples}\label{s:differex}

%The different standards are illustrated in the form of a dialogue.

For examples from topological graph theory and PL topology see \cite[\S2.5]{Sk24} and \cite[\S6]{DS-8}, respectively.

\begin{example}\label{r:ad}
(a) Different standards concerning Remarks \ref{r:frompt}.a.(3), \ref{r:frompt}.b, and \ref{r:fromll}.ab are illustrated by the  discussion of whether the paper \cite{Ad18} is a reliable reference (see explanation in \S\ref{s:intr}) for a proof of the Gr\"unbaum-Kalai-Sarkaria (GKS) conjecture (cf. \cite{Fa20}).\footnote{The most important letters are presented in (d).
Some letters are omitted. In those letters Karim neither stated that according to his reliability standards the paper \cite{Ad18} is a reliable reference for a proof of the GKS conjecture, nor provided a modification that would make the paper \cite{Ad18} such a reliable reference.}
%??? (see though Remark \ref{r:fromll}.a).
%tried to hide (perhaps even from himself) the fact that
%Most of our letters were public, so they are available upon request.
%The letters will be published here only if the judgement of this footnote is questioned, so that a justification by their publication is required. 

An attempt to prove GKS is made in \cite[(3) in p. 7]{Ad18}, where the implication
`\cite[Theorem I(1)]{Ad18} $\Rightarrow$ GKS' is attributed to \cite{Ka91}, and it is written {\it `We will give a simpler, self-contained proof of this implication in Section 4.6'}.
This implication is not explicitly proved

$\bullet$ in \cite{Ka91} (it is even not explicitly stated there) because the statement of
\cite[Theorem I(1)]{Ad18} is not present in \cite{Ka91} (even as a conjecture), and the paper \cite{Ad18} does not mention that \cite[Theorem I(1)]{Ad18} was not known in 1991 (even as a conjecture);

$\bullet$ in \cite[\S4.6]{Ad18} because \cite[\S4.6]{Ad18} does not explicitly mention this implication.

% (see also (b)).

\smallskip
(b) In his letter of 22.01.2019 (but not in \cite{Ad18}, see Remark \ref{r:fromll}.a) Karim Adiprasito writes {\it `...it [AS: GKS] follows from Theorem I and Corollary 4.8. (I also remark that after that corollary).'}
This is wrong because \cite[Corollary 4.8]{Ad18} is not a rigorous mathematical statement.
Indeed, the assumption of \cite[Corollary 4.8]{Ad18} is `\emph{If $\Sigma$ is a $2k$-dimensional rational sphere in $\R^{2k+1}$ that satisfies the hard Lefschetz property}'.
However, the definition of `\emph{the hard Lefschetz property}' is neither presented nor referred-to in \cite[\S1-\S4]{Ad18}.\footnote{Also, an object cannot satisfy a theorem (as in arXiv version 1 of \cite[Corollary 4.8]{Ad18} to which Karim's letter refers; I could not find any statement named `hard Lefschetz theorem' in \cite[\S1-\S4]{Ad18}.}
%Observe that `\emph{the hard Lefschetz property}' is not defined in \cite[Corollary 4.8]{Ad18}, or at least the definition is hard to find in the preceding 23 pages.
Observe that in \cite{Ad18} it is not written that GKS {\it `follows from Theorem I and Corollary 4.8'.}
%Corollary 4.8 of \cite{Ad18} might be something close to this implication, but is not such an implication.

\smallskip
(c) In spite of (a,b) it is still possible that \cite{Ka91} or \cite[\S4.6]{Ad18} contain something very helpful for writing a reliable proof of the implication `\cite[Theorem I(1)]{Ad18} $\Rightarrow$ GKS' (and that the implication should be attributed to \cite{Ka91}).
A reliable reference would explicitly prove the implication (and explain why the implication is attributed to \cite{Ka91}, with references to particular places in the 25-page paper \cite{Ka91}).

In January 2019, when I found the gap described in (a) in version 1 of \cite{Ad18}, I thought the gap is very minor.
As of March 2022, this is a 3-year gap and a half-a-page gap in the sense of Remark \ref{r:fromll}.b.
But I do not assert that the gap cannot be filled.

\smallskip

\small
(d) {\it (AS to KA, 16.01.2019)}

Dear Karim,

Thank you for writing an interesting paper  [v1 of  \cite{Ad18}].  
I am interested in learning your proof of the Gr\"unbaum-Kalai-Sarkaria conjecture, 
but I could not find it in \S1 or \S4. 
If I am missing something, could you let me know in which page(s) the proof is written? 
If not, could you update your paper adding a head  `proof of the Gr\"unbaum-Kalai-Sarkaria conjecture' and such a proof under this head?    

Best, Arkadiy. 

{\it (KA to AS, 22.01.2019; I am grateful to Karim for allowing me to publish this letter)}

Hi Arkadiy,

this is a corollary of Theorem I, and it is somewhat well-known to be (Kalai showed this first in his paper \cite{Ka91}). I did not list all corollaries of that theorem explicitly, because there are too many and the derivations assuming my main theorem are done better by others before me, but the Grunbaum is actually derived explicitly, as it follows from Theorem I and Corollary 4.8. (I also remark that after that corollary).

Best, Karim

{\it (AS to KA, 5.04.2020)}

Dear Karim,

Thank you for your renewed interest in my critical remarks on arXiv:1812.10454.
I sent you all those remarks in January, 2019.
I'm afraid they are not taken into account in \cite{Ad18}.
The proof of GKS seems to be contained in (3) in p. 7, where you attribute the implication
`Theorem I(1) $\Rightarrow$ GKS' to \cite{Ka91} and write
`We will give a simpler, self-contained proof of this implication in Section 4.6'.
Your paper suggests that Theorem I(1) was not known in 1991, even as a conjecture.
Thus implication `Theorem I(1) $\Rightarrow$ GKS' could not be proved in \cite{Ka91}.
Section 4.6 contains no head `Proof of the implication Theorem I(1) $\Rightarrow$ GKS'.
Section 4.6 has only one formally stated result, Corollary 4.8, which is not the implication `Theorem I(1) $\Rightarrow$ GKS'.

It would be nice if you could update arXiv version adding a proof of the implication `Theorem I(1) $\Rightarrow$ GKS' under the name `Proof of the implication `Theorem I(1) $\Rightarrow$ GKS''.
If you send me a project of such an update, I am willing to read it and tell you if I have further questions.

Best, A.
\normalsize
\end{example}

\begin{example}\label{r:fv21}
Different standards concerning Remark \ref{r:steps}.a are illustrated by the following example.
The paper \cite{FV21} refers to \cite{MW16} without indicating unreliability of the proof in \cite{MW16} justified in \cite{Sk17o}, \cite[\S5]{Sk17-3} (and without referring to the different proof of \cite{Sk17-3}).
A discussion with colleagues (or an internet search) would yield \cite{Sk17o, Sk17-3} before arXiv submission of \cite{FV21}.
%https://www.youtube.com/watch?v=3DZ0ETl2PBw&list=PLU5TmElGttj0o7V4oQnRwyx7P5MMcPYh5&index=15  minute 33
The authors of \cite{FV21} were informed about \cite{Sk17o, Sk17-3} and the criticism of \cite{MW16} in June 2021 at Moscow Conference on Combinatorics and Applications.
As of February, 2022, arXiv version of \cite{FV21} is still not updated.
\end{example}

\begin{example}\label{r:stepsex}
Remark \ref{r:steps}.b is illustrated by the following examples:\footnote{When I saw the arXiv papers \cite{MW16, FK17},
% and journal publication \cite{Cu20},
I informed the authors of my opinion expressed in the bullet points and suggested to update the papers.
When the corresponding updates will appear, I would be glad to remove these examples.}

$\bullet$ because of the way my name and work is mentioned in \cite{MW16}, I find it misleading that \cite{MW16} does not mention that a pre-arXiv version of \cite{MW16} was sent by the authors to me, I  liked the idea of proof and had important specific criticism on its realization, but \cite{MW16} was put to arXiv without fulfilling the recommendations \ref{r:steps}.b.(2),(3); cf. \cite{Sk17o}.\footnote{\label{f:help} The papers \cite{Sk17-2, Sk17o} are not results of a research competitive to \cite{MW16}, but are results of several-years attempts to help the authors recovering their gap.
%---------- Forwarded message ---------
%From: arkadiy skopenkov <askopenkov@gmail.com>
%Date: Sun, Jan 15, 2017 at 2:02 PM
%Subject: Re: Eliminating Higher-Multiplicity Intersections
%To: Sergey Avvakumov <s.avvakumov@gmail.com>, Uli Wagner <uli@ist.ac.at>, Isaac Mabillard <imabillard@gmail.com>
See the following letter of A. Skopenkov to S. Avvakumov, U. Wagner, I. Mabillard, from January 15, 2017.
\newline
{\it When I wrote a remark to part III [added in 2020: to \cite{AMSW}], I realized that the idea of recovering the gap in part II [added in 2020: to \cite{MW16}] by smoothing, which I suggested early in 2015, works:
\newline
The smooth version of the Local Disjunction Theorem 1.16 [added in 2020: from version 2 of \cite{AMSW}] is correct.
Moreover, the smooth version implies (by approximation) the
piecewise-smooth and then the PL version.
The same holds for the non-injective version of the metastable Local
Disjunction Lemma \cite[Lemma 10]{MW16}.
Thus by proving the smooth non-injective version of the metastable
Local Disjunction Lemma
(using vector bundles), both the gap (of inaccurate use of PL block
bundles) in \cite{MW16}
can be recovered, and the proof can be shortened.}}

$\bullet$ because of the way M. Skopenkov's work is mentioned in \cite{FK17}, I find it misleading that \cite{FK17} does not mention that pre-arXiv versions of \cite{FK17} were sent by the authors to M. Skopenkov, who liked the idea of proof and sent the authors important specific criticism on its realization, but \cite{FK17}
was put to arXiv without (3); the authors kindly sent some updates to M. Skopenkov, he answered with more critical remarks and encouragement, and the arXiv updates of \cite{FK17} still do not mention that M. Skopenkov did not confirm that his criticism is reasonably resolved.\footnote{I am grateful to Rado Fulek for many fruitful discussions in spite of our differences on this point.}

%the authors sending arXiv versions to M. Skopenkov
%(to check that he approves the way his , if they have been).

%$\bullet$ see the Zentralblatt review on \cite{Cu20}.

$\bullet$ mathematicians would not find it misleading that \cite{Ad18} does not mention that I had some specific criticism on the argument (their opinion might or might not change if they have learned of the criticism, see Example \ref{r:ad}).
%it is sufficient that Karim removed upon my request my name from acknowledgements in version 3 or later.
\end{example}

\begin{example}\label{e:pt19}
Different standards concerning Remark \ref{r:prework} are illustrated by the following example.
The crucial reference to \cite[before Theorem 3.2]{Kr00} is missing in the paper \cite{PT19}.
The only essential `new' point of the general results \cite[\S1.2]{PT19} (as explained in  \cite[Remark 5.d.2]{PT19}) is the calculation of Johnson’s obstruction.
This `new' point is essentially done in \cite[before Theorem 3.2]{Kr00} (as explained in \cite[Remark 2.5.2b]{KS21e}).

\small
Less importantly, the references \cite{FK19, KS21, KS21e} are missing from the paper \cite{PT19}.
Of these,

$\bullet$ paper \cite{FK19} contains a graph-theoretical analogue of \cite{PT19}; the corresponding misleading remarks of \cite{PT19} are criticized in \cite[footnote 3]{Sk24},

$\bullet$ papers \cite{KS21, KS21e} are simplified expositions of \cite{PT19} citing \cite{Kr00, FK19, PT19}.

\emph{Comments.}
These misleading omissions are deliberate.
Indeed, the papers \cite{KS21, KS21e} were sent in 2021 before their arXiv submissions to the authors of \cite{PT19}, with whom we earlier discussed the relevance of \cite{PT19} to the content of \cite{KS21, KS21e}.

See also Math Review of \cite{PT19}.
\end{example}

\normalsize

\section{Reviewing the peer review system}\label{s:peer}

%https://www.ncbi.nlm.nih.gov/pmc/articles/PMC3023113/

\hfill{\it All professions are conspiracies against the laity.}

\hfill{B. Shaw, The Doctor's Dilemma}

\begin{remark}[Why high quality of peer review is vital]\label{r:qual}
Publication in a peer reviewed scientific journal is the main test for recognizing a result as reliable (see  \S\ref{s:intr}).
Thus reliability standards are mostly set up by such journals.
Besides, publications in such journals are used for jobs and grants distribution.
Consequently, such journals practically rule the mathematical world.
%Hence for the existence of a peer reviewed journal it is vital to maintain the following principles
%of Remark \ref{r:prin}.

Sometimes `official' peer review standards differ from practice.
See Examples \ref{r:review}, \ref{r:aks-ijm}, and \S\ref{s:isrask17}, where this is easy to see even to an outsider
who does not has a vast experience of reading and/or writing referee reports.
%For unreliability of the peer review system at an industrial level
See also \cite{CLM, Ab21}.
Consequently, it is vital to recognize and correct such situations, see Remark \ref{r:sugg}.
%, and nullify their influence.

Dissemination of research is not a reason for the existence of a peer reviewed journal
(because research
%papers
can be disseminated by the authors via the internet).
However,
%mathematicians (including computer
scientists might want to attract attention to some results,
but have lack of time and energy to follow the principles of Remark \ref{r:prin}.
Then it is important that the web page of a journal (or a publication) announces that the journal is not a peer reviewed journal, and so should not be used as a reliable reference and/or for jobs or grants evaluation.
Publications in arXiv and computer science conferences are examples of this kind \cite{Fo}.

{\it Without all that a journal degenerates to an instrument of redistribution of jobs and grants in a way obstructing the progress of science.}
\end{remark}

\begin{remark}[Some principles of scientific discussion]\label{r:prin}
(a) {\bf Scientific truth is established by justification and reasoning,} not by authority, majority or administrative pressure.
This is what makes science different from other respectable activities of thought (like religion, intelligence service, etc.).
Thus criticism and responsibility for criticism are vital for scientific discussion.
Without them reliability standards could not be kept.
Suppression of criticism or irresponsible criticism contribute to the degeneration of an activity to a non-science.
Therefore suppression of criticism or irresponsible criticism must be identified as such,
and their negative influence on science must be overcome.

In particular, the following {\it impartiality principle} should be maintained:
decisions ought to be based on objective criteria, rather than on bias, prejudice, or preferring to benefit one person over another for improper reasons.

Since technically the scientific truth is (or should be) established by publications in peer reviewed journals
(see Remark \ref{r:qual}), the same principles should be maintained for editorial decisions in peer reviewed journals.

However, there is nothing wrong if scientists/referees/editors present unjustified  statements, when they  explicitly write that these statements do not affect practical issues (like editorial decision), and are presented on a `take it or leave it' basis.

An example of a competent justification that a paper does not meet the high standards of a particular journal (Fund. Math.) is \cite{Rep} (publication of this report is allowed by the Editor H. Toru\'nczyk).
See examples of unjustified / wrong statements / judgments in referee reports and surveys (including unjustified statements that proofs in published papers are incorrect) in Examples \ref{r:ps-dcg}.b, \ref{r:aks-ijm}, \ref{r:revrevi} and \cite[Remarks 5.4, 5.5, 5.6, 5.9]{Sk16}.

%\begin{remark}[Handling of \cite{KS20} in Fund. Math]\label{r:ksfm}
%The referee report \cite{Rep} (whose publication is allowed by the Editor) is an example of a
%competent justification that the paper does not meet the very high standards of a particular journal.
%Unfortunately, such competent justifications are so rare nowadays.
%The paper was rejected.
%\end{remark}

(b) {\bf Confidential decisions should match public discussions.}
Every substantial scientific argument affecting a confidential decision should be publicly available.\footnote{Here a `decision' is a decision by an editorial board, dissertaion/diploma committee, etc.
The {\it quantity} of persons publicly supporting some point of view need not be matched by a confidential decision.}
An example on dissertation defence in Russia: if there are no public objections to a dissertation, but a dissertation committee votes against the dissertation, then the committee is dissolved.
(I hope there are analogous rules or traditions in other countries.)

This is required for maintaining (a).
In particular,

$\bullet$ criticism and different opinions should not be suppressed by administrative means
(e.g. misusing the anonymous review system by writing biased negative reports on papers whose authors propose  criticism and opinions different from the reviewer's);

$\bullet$ references to criticism and to different opinions should be cited but not suppressed
(even if one disagrees with them).

{\it I disapprove of what you say, but I will defend to the death your right to say it.}

%I detest what you write, but I would give my life to make it possible for you to continue to write.
%Evelyn Beatrice

(The Friends of Voltaire, E. B. Hall)

See \cite[Remark 5.1.b]{Sk16}; e.g. compare \cite[\S1]{Sk16} citing \cite{BZ16}, and \cite{Sh18} citing \cite{Sk16} to \cite[\S1]{BZ16}, \cite{BS17, FS20} not citing \cite{Sk16}.

The anonymous review system gives a convenient frame for suppressing criticism, for irresponsible criticism, and for promotion of unreasonable opinions which do not withstand open discussion.
See Examples \ref{r:skrms}, \ref{r:skajm}, in \S\ref{s:isrask17}, and in  \cite[Remark 5.3 and \S5.4]{Sk16}.

(c) {\bf Decisions should be made carefully and critically studying the texts in question,} not blindly believing (charming or influential) people we know, cf. \cite{Gr}.
The texts in question are

$\bullet$ referee reports and authors' responses, for editorial decisions,

$\bullet$ research and survey papers, for distribution of credits.

This is required for maintaining (a,b).

(d) See examples explicitly showing that the principles (a,b,c) are maintained (sometimes in spite of attempts to violate them) in Examples \ref{r:skrms}, \ref{r:csmmj}, \ref{r:skajm}, \ref{r:ps-ta}.
These examples represent high standard of impartiality and freedom of scientific discussions maintained by the Editors.
\end{remark}

\begin{remark}[Mistakes in editorial decisions]\label{r:sugg}
(a) With all our concern for reliability, mistakes are unavoidable.
Thus a characteristic of a peer reviewed journal is that the Editors recognize and correct their mistakes.
%This is vital to show that editorial decisions are made with reasonable amount of care (see \S\ref{s:differ}
%and Remark \ref{r:prin});
Since most of the Editors' work is not immediately visible to a reader, failure to recognize and correct their mistakes even in a few cases indicates that the degeneration described at the end of Remark \ref{r:qual}  takes place.

%Unless explicitly indicated otherwise, a medical drug should work properly in overwhelming majority of cases

If a mistake involves {\it acceptance} of a paper containing unreliable result/proof, then the correction involves publication of errata in the same journal (whether by authors, or by Editors, or by a third party).

(b) If a mistake involves {\it rejection} of a paper because of an unreasonable referee report, I suggest the following procedure for keeping (or restoring) the journal's status of a peer reviewed journal.
The rejection decision should be canceled, and consideration of the paper should be renewed.
The paper should be accepted or rejected only upon a referee report (justifying its conclusions) sent to the author, not upon the Editors' opinion without sending a referee report.\footnote{Sometimes the initial submissions are rejected upon the Editors' opinion that the paper is not suitable for the journal.
This policy has its advantages (referees do not waste their time on written explanations that the paper is not suitable; authors can soon send the paper to a different journal) and drawbacks (decision without a report is potentially unreasonable or even corrupt). It is outside our purposes to discuss this policy here, we only state that such a policy is clearly inappropriate after rejection of a paper because of an unreasonable referee report.
If one report is unreasonable, and other reports are hidden, a scientist has to suppose that the other reports are also unreasonable.}
Indeed, even a great expert can either make a mistake, or be biased, or be emotionally dependent on a less qualified/honest person (see Remark \ref{r:prin}.c and \cite[\S5]{Sk16}).
So without sending a referee report to the author, the decision is way too unreliable for the status of a peer reviewed journal (described in Remark \ref{r:qual}).

The conclusion that a referee report is unreasonable has to be ignored if it is not justified by specific criticism of specific phrases of the report.
The same holds for any other conclusion in the discussion described below.
If the conclusion is justified by specific criticism, then it should be sent to the referee, so that he/she could provide a revised report (taking into account the criticism), or explain why he/she disagrees with the justification, or refuse to reply (in the last case the conclusion is recognized to be correct).
The referee's letter is sent to a person who justified that the referee's report is unreasonable (usually to the authors), and so on.
The Editors can invite alternative referee for the paper or for the particular questions discussed.
When the discussion converges, the Editors make their new rejection/acceptance decision.
For an effective means to make the discussion responsible (and hence short) see Remark \ref{r:transp}.b.
If the Editors feel lack of time and energy to moderate such a discussion, {\it `it is important that the web page of a journal announces that the journal is not a peer reviewed journal, and so should not be used as a reliable reference and/or for jobs or grants evaluation'}, see Remark \ref{r:qual}.\footnote{The main purpose of this  procedure is keeping (or restoring) the journal's status as a peer reviewed journal.
If the authors wish that their paper (rejected because of an unreasonable referee report) should be reconsidered in the same journal, then this purpose is achieved together with the reconsideration.
If the authors wish to submit their paper to a different journal, the above procedure still makes sense (because of its main purpose).
Then  the formal result of such a procedure is `honorary rejection' or `honorary acceptance' of the paper.
The journal publishes the information on the result but not the paper, in order to comply with the copyright. The result is `honorary reconsideration after revision' if the authors do not submit a revised version after a reasonable amount of time specified by the Editors.}
\end{remark}

\begin{remark}[Transparency helps]\label{r:transp}
(a) The transparent anonymous peer review policy \cite{Skt} (perhaps introduced gradually and partly)
would be helpful for maintaining the principles of scientific discussion (described in Remark \ref{r:prin}).
Transparency allows to diminish the gap between official information and practice, see \S\ref{s:peerex} and \S\ref{s:isrask17}.

(b) The transparent anonymous peer review policy \cite{Skt} allows not to waste our
%Editors', referees', and authors'
valuable time on replying to incompetent texts.
Those who write incompetent texts usually are competent enough not to allow publication of these texts on the internet, and such a policy involves ignoring texts not put on the internet.

(c) When a discussion is an argument rather than a collaboration, we strongly need it to be responsible.
In such situations we usually do not have enough time to discuss premature ideas, whose invalidity becomes clear when their publication (or a mental experiment of publication) is suggested.
(Such ideas usually, though not always, are euphemisms for `I have more administrative power than you'.)
Then it is useful to have a transparent discussion, in which letters are (allowed to be) published on the internet.
If a participant of such a transparent discussion receives a letter not stated to be public,
then he/she deletes  it unread (to avoid confusion).
If a part of such a public discussion becomes obsolete, participants can delete that part
(only) by mutual consent.

(d) {\it `Public shaming is the only thing that can really work against groupthink. To spread the word, please LIKE this post, REPOST it, here on WP, on FB, on G+, forward it by email, or do wherever you think appropriate.'} \cite{Pan}
\end{remark}

\section{Reviewing the peer review system: examples}\label{s:peerex}

%I present my experience as an author.
{\bf Conventions for \S\S \ref{s:peerex}--\ref{s:letters}.}
In the sources (i.e. in the letters and in citations from referee reports) the references are updated.
In spite of that, the references are given to the very version of a paper cited in the source.
Otherwise the sources are not changed (the grammar is not corrected).
Quotations from referee reports are given in italics and in quotation marks.
We quote all parts of the reports related to the rejection recommendation.
We omit e.g. quotation of referee's description of current situation in the area before the appearance of the paper under review.
For those referee reports that are both incompetent and play a crucial role in the Editors' decision,
the complete texts are available upon request (in particular, to confirm that we did not omit
any critical remarks justifying the rejection recommendation).
If an unreasonable report did not play a crucial role for the Editors' decision, there is no need to justify my criticism by publishing my reply to the report.
I omit cases before 2017 and those cases which did not involve a non-trivial discussion (i.e., those cases in which a paper was accepted after a reasonable amount of changes required by the referee).

%[Official information and practice]

\begin{example}\label{r:review}
`Official' peer review standards include reasonable checking of novelty of results.
However, I do not know any peer reviewed journal which requires submitted papers to be publicly available on arXiv.
If a paper is not publicly available on arXiv, then I consider novelty of its results as not being reasonably checked, see Remark \ref{r:steps}.
%I do not know any opposite public statement by a journal (or a journal editor).
%clarity and novelty
Observe that checking of novelty in the form putting a text on arXiv is necessary (although not sufficient)
for a paper to be awarded with a {\it scientific prize} of the Moscow Mathematical Conference of High School Students, see \url{https://mccme.ru/mmks/index_eng.htm}.
\end{example}

\begin{example}[the report on handling of (a previous version of) \cite{Sk16} in Russian Math. Surveys]\label{r:skrms}
One of the referee reports misused the anonymous peer review system to promote an unreasonable opinion which does not withstand open discussion, see \cite[\S5]{Sk16}.
{\it I am grateful} to Editor S. Shlosman for {\it publicly criticizing} in \cite{Sh18} the survey \cite{Sk16} because this allows to see that the criticism exists and is incompetent (see \cite[Remark 5.6]{Sk16}).
{\it I am grateful} to the Editors for publication of \cite{Sk16} without correcting the passages criticized by  S. Shlosman, and publishing S. Shlosman's criticism in the same volume \cite{Sh18}.
\end{example}

\begin{example}[the report on handling of \cite{CS16} in Algebr. and Geom. Topol.]\label{r:csagt}
To appear.
\end{example}

\begin{example}[the report on handling of \cite{CS16} in Moscow Math. J.; A. Skopenkov]\label{r:csmmj}
One of the referees (Alexey Zhubr) kindly disclosed his identity to facilitate the discussion of the paper.
He presented important criticism and justly asked for a major revision.

In response to his report made after such a revision
%he criticized the paper for unclear exposition and provided a recommendation (however reserved)
%to publish this paper.
the authors wrote to the Editors,

`The referee of our paper Alexey Zhubr kindly forwarded to us his report of 19.09.2019.
We are grateful to him for a thorough reading of our paper, many specific suggestions he made in earlier reports, and a recommendation (however reserved) to publish this paper.

The critical conclusions of the latest report are not justified by references to the paper and suggestions what could be done in a more clear way.
So we can only state that we disagree with that conclusions, that sections 3 and 4 are mostly hard to read because of complexity of the matter (not because of poor exposition), and that we invested several years to write this paper in a clear way.
In our opinion, a referee should be requested to justify his critical conclusions
only if the Editors think the conclusions could affect the acceptance of a paper.
Otherwise the matter could be dropped.'

Another referee presented minor but valuable critical remarks.
The paper was accepted without further revision.
\end{example}

%(d) {\it A. Skopenkov's review of handling of \cite{ST17} in Discr. Comp. Geom.}
%The report(s) arrived reasonably soon, they contained a minor criticism and asked for a minor revision.

%\footnote{See  \url{https://scirev.org/journal/arnold-mathematical-journal/}

\begin{example}[the report on handling of \cite{Sk18} in Arnold Math. Journal]\label{r:skajm}
There were two initial reports on the paper.
The first was positive, it contained important criticism and asked for a major revision.
The second was negative, it misused the anonymous peer review system to promote an unreasonable opinion which does not withstand open discussion, see \cite[\S5]{Sk16}.
%https://arxiv.org/abs/1605.05141).
In my reply to the Editors, I justified the latter judgement by considering the referee's comments one by one.
The Editors suggested a major revision (even before receiving my criticism of the second report).
I do not know whether the Editors sent my criticism to the second referee or not.
I received no reply to my criticism from the second referee.
The paper was accepted after several iterations of work upon suggestions of the first referee.
%Since the unreasonable report did not play a crucial role for the Editors' decision, there is no need to
%justify my criticism by publishing my reply to the report.
\end{example}

\begin{example}[the report on handling of \cite{PS20} in Discr. Comp. Geom.]\label{r:ps-dcg}
Our review is presented in the form of our letter of 19.12.2021 to the Editors J. Pach and K. Clarkson.

\smallskip
(a) {\it Dear Janos, Dear Kenneth,}

Hope you are fine and healthy.

Please find after our signatures [added later: in (b)] our report to review \#2 on our paper
\cite{PS20} recently rejected from DCG.
The report shows that the review is incompetent.

This paper is not among our best papers.
Knowing how busy you are, we are not asking you to reconsider the rejection decision.
However, it is important to recognize incompetent reports as such, and nullify their influence.
Now that papers can be distributed via the internet, the only reason
for the existence of scientific journals is high-quality peer review.
So it would be nice if you could

(1) put the incompetent review to databases (used by other journals) only together with our report.

(2) be more careful with the reviews of the same referee to other papers.

(3) send our report to the referee so that he/she could see that the incompetence is checked
(or could object to our report, if he/she could).

We will not object if you would choose to reconsider the rejection decision.
The other review \#1 is positive. It leaves to the Editorial Board the judgment
whether the paper is sufficiently related to DCG. This judgment is done by
regular publications in DCG of paper on algebraic and geometric topology
(related to discrete mathematics), e.g. arXiv:1302.2370, arXiv:1307.6444,
arXiv:1512.05164, arXiv:1703.06305, arXiv:1808.08363. Please also observe that
our paper was cited (even before publication) in the paper arXiv:2108.02585
on discrete geometry and computer science.

[A technical paragraph omitted]

\iffalse

A. Skopenkov's reviews of handling his papers are published in \S\ref{s:peerex}, \S\ref{s:isrask17} of \cite{Sk21d'}.
%arXiv:2101.03745
A. Skopenkov plans to add the above description of review \#1.
A report of the review \#2 given below will be added

--- in a complete form (except the paragraph on `minor errors'), if this report contributes to the final Editorial decision;

--- in a significantly shorter form, otherwise.

Reference on the paper
On embeddability of joins and their `factors'. by
S. PARSA AND A. SKOPENKOV

This article is devoted to the problems of embeddability of polyhedra in Euclidean
spaces. Thus one of the main results of the article establishes the equivalence
of the embeddability of polyhedra in a Euclidean space to the embeddability of
the Union of its three cones in the product of this space on the plane The main
results are clear, convincingly proven and have the necessary novelty. In General,
this is a good article on geometric topology. However, its relation to discrete or
computational geometry is problematic. Whether it should be published in your magazine is up to the editorial Board.

\fi

Best Regards, Salman and Arkadiy

%Reviewer #2: In my opinion  the paper should not be considered for publication in DCG since it does not meat
%the standards. The paper presents a proof of a particular instance of a result by Melikhow—Schepin.
%Furthermore, the paper looks like a "piece/subsection taken out of an actual research paper".
%The result presented in this paper  (including the proof techniques) is not of  broad interest for
%the readers of DCG and will not attract desired attention. For these reason I suggest rejection of the paper.

\smallskip
(b) {\it A report to the review \#2:}

The first sentence of our paper is

`{\it We present short and clear proofs of Theorems 1 and 3.b which first appeared in the unpublished paper
\cite[(iv)$\Rightarrow$(i) of Corollary 4.4, Theorem 4.5]{MS06}.}'

The second sentence of the review {\it `The paper presents a proof of a particular instance of a result by Melikhow-Schepin'} repeats the above first sentence of the paper in a misleading form.
(The `short and clear proofs' [of the current paper] and `unpublished' [paper of Melikhov-Schepin] are missing.)

[a paragraph on `minor errors' in another paper is omitted]

%The proofs of the results from the Melikhov-Schepin paper just cited contain `minor errors'
%recognized by Melikhov in a discussion with the first author.
%In spite of that, we do not intend to have a priority argument or to mention the errors in our paper.
%In a future paper, whose draft was sent by Melikhov to the first author, Melikhov will fix these
%and provide better proofs of the unpublished paper.

The third sentence `{\it Furthermore, the paper looks like a ``piece/subsection taken out of an actual research paper''}' instead of justifying the judgment of the review, gives the paper a negatively sounding name (`piece/subsection') insinuating that the paper is not a research paper.
However, the third sentence does not explicitly state that the paper is not a research paper (because it clearly is), does not explain why the paper is not a research paper, does not describe what is `piece/subsection', and
does not explain why the paper is such.
The reviewer did not refer to the actual paper of which our paper is a subsection.
Every short (or long) paper could be in principle a part of a larger paper.
Thus the third sentence is a logical fallacy.

The remaining three sentences of the review are (partly repeated) judgements based on the above-discussed second and third sentences.
Thus the review attempts to justify that the paper does not meet the standards of this journal,
but fails to do so.
\end{example}

%not published

\begin{example}[the report on handling of \cite{PS20} in Topol. Appl.]\label{r:ps-ta}
The initial report was incompetent.
This was shown in our reply to the report, and this was not questioned (cf. Example \ref{r:ps-dcg}).
The Editors considered different opinions on the paper.
The paper was accepted after work upon suggestions of an alternative referee.
\end{example}

\section{Violation of important scientific principles}\label{s:isrask17}

\begin{example}[introduction]\label{r:viol}
\textbf{(a)} This section (and \cite[\S6]{DS-8}) justifies that principles of scientific discussion (described in Remark \ref{r:prin} and in \cite[Example 6.9]{DS-8}) are violated in
Bull. of the London Math. Soc.,
Fund. Math.,
Intern. Math. Res. Notes,
Israel J. Math., and
Proc. of the Amer. Math. Soc.

The papers \cite{Sk17-2}, \cite{AKS-1}, \cite{AKS-2}, \cite{KS21e}, \cite{DS-3},
\cite{DS-5} (more precisely,  their shortened versions as explained in Examples listed below) were rejected from those journals (except Bull. of the London Math. Soc.), based on incompetent reports.
The incompetence is justified in Examples \ref{r:revrevi}--\ref{r:dsirmn5}
%, \ref{r:aks-ijm}, \ref{r:akspamsd}, \ref{r:ksfm},  \ref{r:dsirmn}, and \ref{r:dsirmn5}
(see also our letters in \S\ref{s:letters}).
The justifications was sent to the Editors.
The validity of justifications is partly confirmed by the lack of replies to them.
%The reaction of Israel J. Math., whenever available, will be published here.
This lack of replies
%reaction will hopefully contribute to
shows that taking editorial decisions on the basis of incompetent referee reports is a consistent policy of the journals.
%See \S\ref{s:peer} and \S\ref{s:isrask17}.
Although the justification is not questioned by the Editors, the rejection decision is not withdrawn (see also footnotes \ref{f:ijm} and \ref{f:dsirmn}).
This violates the principles of Remark \ref{r:prin}.ab.
\footnote{No public reply to our letters of Remarks  \ref{r:ijm1}--\ref{r:dsirmn-l} are available by July, 2025.
Therefore these corresponding examples are gradually moved from  \S\ref{s:peerex} to \S\ref{s:isrask17}.
Public replies, if available, will be presented in an update of this paper.}

The papers \cite{Sko, Skm} and Example \ref{e:moebius} present examples of suppression of criticism by administrative means, and of the difference between official information and practice, violating these principles.
So

$\bullet$ Bull. of the London Math. Soc.,
Fund. Math.,
Intern. Math. Res. Notes,
Israel J. Math., and
Proc. of the Amer. Math. Soc.
are not peer reviewed {\it scientific} journals,

$\bullet$ Department of Differential Geometry of Faculty of Mechanics and Mathematics of Moscow State University is
%currently
not a {\it scientific} department, and

$\bullet$ the M\"obius contest is not a {\it scientific} contest.

Instead, these organizations are instruments of redistribution of jobs and grants in a way obstructing the progress of science.
%??? (cf. (c)).

It is important for the community to know that.
Also, this would hopefully be useful for the Editors of those journals, either

$\bullet$ to restore the status of a peer reviewed scientific journal (see Remark \ref{r:sugg}), or

$\bullet$ to publicly explain why they disagree with specific parts of Remarks \ref{r:qual}, \ref{r:prin}, \ref{r:sugg}, or

$\bullet$ to announce at the web page of the journal that {\it `the journal is not a peer reviewed journal, and so should not be used as a reliable reference and/or for jobs or grants distribution'} (see Remark \ref{r:qual}).

For analogous reason this would hopefully be useful for the chair and members of the department, and for the   president and members of the contest.

\smallskip
\textbf{(b)} {\it `Violation of important scientific principles ... is a consistent policy of an influential group in the community of topological combinatorics'} \cite[Remark 5.3]{Sk16} (see justification in \cite[Remark 5.3]{Sk16}).
Material which is reliable, openly unopposed, and referring to the different (however misleading) point of view is not published, while material which is misleading, publicly criticized, and failing to refer to public criticism, is published, see \cite[Remark 5.3]{Sk16}.
This jeopardizes the peer review system.

Examples \ref{r:skrms}, \ref{r:skajm}, in \S\ref{s:isrask17}, in \cite[\S5]{Sk16} show that this group (or mathematicians blindly believing private judgements or publications by this group) attempts to misuse journals %(believed to be peer reviewed journals)
as instruments of redistribution of jobs and grants in a way obstructing the progress of science.
The examples are different: Russian Math. Surveys and Arnold J. Math. are not such instruments, some journals have a chance to correct mistakes exhibited to Editors (see Remark \ref{r:sugg}.a), while Israel J. Math. and Proc. of the Amer. Math. Soc. remain such instruments even after mistakes are exhibited to the Editors.

%(c) !!!
\end{example}

Recall the conventions from the beginning of \S\ref{s:peerex}.

\begin{example}[the report on handling of \cite{Sk17-2} in Israel J. Math.; 29.01.2021]\label{r:revrevi}
The paper was rejected on the basis of a referee report shown to be incompetent below.
%(cf. footnote \ref{f:ijm} and Example \ref{r:ijm1}).
The referee's specific critical remarks are either

$\bullet$ attempts to misuse the anonymous peer review system to promote an opinion which does not stand an open discussion (and thus violating principles of scientific discussions recalled in Remark \ref{r:prin}.b), see quotations in (A,E) below, or

$\bullet$ misleading judgements based on incomplete quotations of the manuscript, see quotation in (D) below, or

$\bullet$ statements and judgements irrelevant to the rejection recommendation, see quotations in (C,F) below.

Most important replies are (A)-(F); more technical replies are (F1)-(F5) and (G).\footnote{\label{f:ijm} The rejection decision on \cite{Sk17-2} was confirmed, referring to other experts' opinion, but without presenting any new referee report; no report was sent upon my request of Remark \ref{r:ijm1}.b.
A private Editor's reply to my letter of Remark \ref{r:akspamsl}.a did not demonstrate understanding of this problem, and did not suggest or accept any way of resolving the problem.
\newline
Additionally, the impartiality principle of Remark \ref{r:prin}.b is violated for \cite{Sk17-2} because a similar paper \cite{AMSW} containing a {\it weaker} result was published in the Israel J. Math.}

\smallskip
{\bf (A)} `{\it 1. The current manuscript does not give any reasons why the author finds the
proof in \cite{MW16} unsatisfactory. Instead the reader is pointed to reference
\cite{Sk17o}, an arXiv preprint that quotes emails between the author and Wagner. In
\cite{Sk17o} the author gives a list of reasons why the proof of \cite{MW16} should be
regarded as incomplete. Even if the treatment in \cite{Sk17o} was satisfactory and
well-reasoned (it is far from that), since it is in a different manuscript that
presumably is not undergoing peer-review, the reader is simply asked to take the
author's word for the incompleteness of \cite{MW16}. This is unacceptable...}'

%The reader is referred to [Sk17o, Remarks 3 and 4] for {\bf detailed description of problems with [MW16]} and for a reply to that description by one of the authors of [MW16].
%$\bullet$ The paper [Sk17o] presents not only criticism of [MW16], but also the reply of an author of [MW16] to this criticism.

The referee does not justify the judgement that the treatment in \cite{Sk17o} is far from being satisfactory and well-reasoned (presumably because the treatment is satisfactory and well-reasoned).
This is too strong a judgement to be used without justification for rejection recommendation.
The more so because an author of \cite{MW16} wrote `I agree with many of the specific criticisms ... I think these are helpful and valid criticisms, and we will address them in the next revision' \cite[Remark 4]{Sk17o}, and then no next revision appeared in 3 years.
Cf.  arXiv:2101.03745v1, Remark 1.3.b.

{\it `The reader is simply asked to take the author's word for the incompleteness of \cite{MW16}'}

is wrong because

{\it `In \cite{Sk17o} the author gives a list of reasons why the proof of \cite{MW16} should be regarded as incomplete.'}

In any case, these remarks of referee cannot justify the rejection recommendation because of (B) below.

I did not include the relevant part of \cite{Sk17o} into the submission because of (B) and (F3) below.

%I am resubmitting the paper keeping [Sk17o] as an external reference and asking the referee to present his specific critical remarks to [Sk17o, Remark 3] if he/she has any.
%(I do this on the assumption that Editors will find that the relevant part of [Sk17o] would not be interesting to most readers.)
%However, I am willing to fulfill any Editors' suggestion on the relevant part of [Sk17o], e.g. to add the relevant part of [Sk17o] to the submission.
%Such a suggestion would guarantee that the increased length of the paper will not be used to justify the rejection decision.

\smallskip
{\bf (B)} The manuscript makes the priority questions irrelevant to its publication by naming Theorem 1.2 Metastable Mabillard-Wagner Theorem.

%The submitted manuscript gives a complete proof of Theorem 1.2.

Remark 1.4 of the manuscript \cite{Sk17-2} explains in mathematical terms (i.e. without an attempt
to distribute credits) the relation of the manuscript to earlier papers (including \cite{MW16}).
This remark is not questioned by the referee.

%In my opinion, journals exist mostly for readers.
A reader is mostly interested in having a reliable reference, not in distribution of credits.
%Of course proper distribution of credits becomes important if there is a priority argument.
%However, there is none in this case:
%The disputed question is whether [MW16] is a reliable reference or not for the result called in the submitted paper {\it Metastable Mabillard-Wagner Theorem.
%So no priority dispute is raised in the submitted paper.
In my opinion, the proof in the manuscript is so short (comparatively to the proof of [MW16]) that it is interesting to a reader independently of his/her judgement on
%distribution of credits (in particular, independently on his/her judgement on
how serious are the gaps in \cite{MW16} described in \cite{Sk17o}.

No part of \cite{Sk17-2} concerns the question whether \cite{MW16} is a reliable reference or not.
Only \S5 of \cite{Sk17-3} concerns this question.
I added of \cite[\S5]{Sk17-3} upon criticism of the referee (see (A)) and I am ready to remove it upon request.
%if the Editors would find it efficient to ask a referee to concentrate on the paper under review, not on [MW16].

%$\bullet$ `In spite of all that I call Theorem 1.2 Metastable Mabillard-Wagner Theorem, in order to concentrate on mathematics and on reliability standards
%for research papers, not on priority question.'
%(Remark 1.4.b of the current manuscript and [Sk17o, p.2].)

\smallskip
{\bf (C)} `{\it 1. If the current manuscript is not supposed to be regarded as an alternative proof but
instead as the first complete proof of the main result, a critical treatment of \cite{MW16} and its perceived gaps is unavoidable.

2. There is a possibility that the current manuscript takes credit away from the authors of \cite{MW16}...}'

These remarks cannot justify the rejection recommendation (and the second sentence is incorrect) because of (B).

\smallskip
{\bf (D)} P. 1 [added in 2022: of the report]. `{\it The proof given in the present manuscript follows roughly the same ideas as in \cite{MW16} but is different:...}

{\it 2. ... The manuscript under review borrows proof ideas (and the statement of the main result) from \cite{MW16}.  For example, in Remark 1.4(b) one reads "Proofs in this paper and in \cite{MW16} are similar...".}'

These judgements are misleading (at least without reference to detailed explanation in  Remark 1.4(b)) and are based on incomplete citation.
Indeed Remark, 1.4(b) reads

`Proofs in this paper and in \cite{MW16} are similar  because they use and extend known methods.
Among these methods are

$\bullet$ for Theorems 1.1 and 3.1: Surgery of Intersection Lemma 2.2 (in \cite{MW16} this method was used without the references presented in the proof of the Surgery of Intersection Lemma 2.2 but with a reference to the `predecessor' paper \cite{Mi61}); ...

The new parts of proofs in this paper and in \cite{MW16} are essential and are different.'

See also \cite[footnote 5]{Sk17-3}:

`...In \cite[\S4.1]{MW16} this surgery of double intersection was used with a reference to the `predecessor' paper \cite{Mi61}, but without the above references...
Lack of the above references in \cite{MW16} may result in exaggerating the overlap of methods of \cite{MW16} and of this paper.'

\smallskip
{\bf (E)} `{\it 2. ... Thus this will have to be handled differently, so that the authors of \cite{MW16} may receive due credit.}

The referee's justification  that the current manuscript does not give due credit to \cite{MW16} is incorrect, see (B) and (D) above.
However, I am flexible in terms of credit distribution (although I am keen on mathematics behind credit distribution).
So I am willing to give more credit to \cite{MW16} upon referee's suggestion of a specific phrase
(if he/she thinks that \cite[Remarks 1.4.ab]{Sk17-2} are not proper or insufficient).
The referee did not state any such suggestion.

\smallskip
{\bf (F)} `{\it 3. The expository quality of this manuscript is low. For example,

--- Several notions are not defined, and the reader is not directed towards a
suitable reference. Examples include "stably parallelizable manifold" or
"$\Sigma_r$-equivariant map," where the notation $\Sigma_r$ is not introduced.

--- The reader is not guided through proofs in a coherent manner. The proof of Lemma 2.2 ends with "So the
lemma is proved analogously to \cite[Theorem 4.5 and appendix A]{HK98},
\cite[Theorem 4.7 and appendix]{CRS}, cf. \cite{Ha63}, \cite[Lemma 4.2]{Ha84},
\cite[proof of Theorem 1.1 in p. 7]{Me17} \cite[The Surgery of Intersection Lemma 2.1]{AMSW}; for detailed descriptions of the above references see \cite[Remark 2.3.a]{AMSW}." This is not helpful exposition.}'

%(Here ? is my typo, it should be [CRS11].)

`{\it --- On p. 5 one reads "By the Hirsch-Smale-Gromov density principle \cite[Corollary in \S1.2.2.A]{Gr86}..." This should be stated explicitly.}'

%(Here ? is the referee's typo, it is \S\  in the paper.)

`{\it --- On p. 3 one reads "The following result implies the existence of a polynomial
algorithm for checking almost $r$-embeddability in $R^d$ for fixed $k, d, r$ such that
$rd \ge (r + 1)k + 3$, cf. \cite[Corollary 5]{MW16}." This should be explained. It
seems the work in preparation referenced in \cite{MW16} never appeared even in preprint form.}'

%(I apologize for a typo: in the paper, and so in the report, it should be [MW16, Corollary 5] not
%[MW15, Corollary 5].)

`{\it --- Abstract and p. 2: "The r-fold analogues of Whitney trick were ?in the air?
since 1960s." What does that mean? The r-fold analogues of the Whitney trick were
formulated and proved by Mabillard and Wagner.}

[The above citation of the report is added later]

The suggestions on exposition in item 3 of the report are reasonable or at least acceptable (except the last one, see (F5) below).
However, these suggestions cannot contribute to the rejection recommendation because

$\bullet$ the referee did not describe any problems with applying well-known argument for proving Lemma 2.2;

$\bullet$ the suggestions are easy to fulfill, see the attached version;

$\bullet$ of (F1)-(F5) below.

In \cite{Sk17-3} I added all the details the referee asked (including a proof of Lemma 2.2).
I am willing to remove them upon request.
% of a referee or the Editors.

%(on the assumption that the details of this well-known argument would not be interesting to a reader).
%However, I am willing to fulfill any Editors' suggestion on this proof (e.g. to add 0.5-1 page of the details of this well-known argument).
%Such a suggestion would guarantee that the increased length of the paper will not be used to justify the rejection decision.

The referee's critical remarks from item 3 should be compared to \cite[Remark 3]{Sk17o}.
The referee's using item 3 for rejection recommendation should be compared to his/her judgement that the criticism
in \cite[Remark 3]{Sk17o} is far from being satisfactory and well-reasoned, see (A).

%\centerline{\bf Other replies to specific remarks by the referee}
%P. 1. `{\it The main result of the manuscript under review is the same as that of [MW16]. }'

\smallskip
\small

Items (F1)-(F5) below reply to the suggestions on exposition in item 3 of the report.

(F1) The notion of a stably parallelizable manifold is a textbook notion.
% well-known to differential topologists.

(F2) The notation $\Sigma_r$ is well-known and/or easy-to-guess to mathematicians working in the area.

(F3) Journals usually have restrictions on lengths of papers.
One of my papers was accepted only on the condition that some details of the argument will be left out and kept only in arXiv version.
Referee's recommendations to give more (or to give less) details are welcome.
Recommendation to give more details can contribute to the rejection recommendation only if the referee justifies that the missing details are important
(which is not so in the current report).
There is nothing wrong that different referees in different journals have opposite recommendations,
but only if such recommendations are not used to justify the rejection decision.
The rejection decision should be based on the mathematical quality of a paper, not on the author's ability to predict the referee's wish to have more (or to have less) details.

(F4) The sentence citing \cite[Corollary 5]{MW16} does not affect the main result.

(F5) Why `the r-fold analogues of Whitney trick were in the air since 1960s'
is explained in the survey \cite{Sk16} cited at the end of the relevant paragraph of \cite{Sk17-2},
see \cite[Remark 3.6]{Sk16}.
In 2017 U. Wagner kindly confirmed that he thinks this remark is proper.

%In the resubmitted version for the sake of non-specialists
%in the area I added the above detailed citation of [Sk16].

\smallskip
{\bf (G)} `{\it 2. ...If, however, there is a serious flaw in [MW16], then the as of now unpublished [MW16] should not appear in its current form.}'

%My English is not sufficient to understand this phrase.
This is unclear because this is not a correct English sentence.
Presumably this phrase refers to an update of [MW16].
Since such an update is not available to math community, the phrase is not relevant to the review.
\normalsize
\end{example}

\begin{example}[the report on handling of \cite{AKS-1} in Israel J. Math.]
%17.03.2022
\label{r:aks-ijm} (R. Karasev and A. Skopenkov)
The paper was rejected on the basis of a referee report shown to be incompetent below.
% (cf. R \ref{r:ijm-ks}).
The referee's specific critical remarks are either

$\bullet$ attempts to misuse the anonymous peer review system to promote an opinion which does not stand an open discussion (and thus violating principles of scientific discussions recalled in Remark \ref{r:prin}.b), see (c,e1,t1,v1) below, or

$\bullet$ wrong or unjustified statements and judgements, see (a,d,e,f,g,e1,f1,j1,k1,l1,o1,p1,q1, r1,s1,t1,u1,v1,w1) below, or

$\bullet$ attempts to create an impression of presenting a justification, without
explicitly calling a justification something which does not stand the test of being called such,\footnote{This could be not an intentional attempt to mislead the Editors, but lack of responsibility, and so lack of understanding the necessity of properly structuring the report, explicitly writing which judgements are justified by which specific remarks.} see (b,g,e1,l1,p1,w1) below.

$\bullet$ statements and judgements irrelevant to the rejection recommendation, see (d,b1,c1,f1,g1, m1,u1,x1,y1) below, or

$\bullet$ valid critical remarks that cannot contribute to the rejection recommendation, see (d1,n1,w1) below.

The most important critical remarks to the report are (a)-(g) below; more technical remarks are (b1)-(y1) below
(the numeration starts with b1 not a1 for technical reasons).

Below the notation and references in the referee report are changed to agree with the notation of \cite{AKS-1},
so as to avoid confusion; otherwise the report is not changed.
Page references are to pages of the report.
In spite of our criticism, in some places we would agree to change phrasing, but the referee does not provide any suggestion (except for suggestions violating the principles of scientific discussions recalled in Remark \ref{r:prin}.b, see e.g. (c) below).

%{\bf Most important replies to specific remarks by the referee.}

\smallskip
{\bf (a)} P. 3. `{\it The main result of the paper, Theorem 1, for certain class of parameters improves the gap (1) as a function of `dimension' and `number of intersections' compared to the gap given in \cite{BFZ}.
This improvement does not give any new insight on \cite[Conjecture 5.5]{BFZ}.
Hence, the level of the results of the paper does not meet the high standards of the Israel Journal of Mathematics
and so, based only on this fact, {\bf I suggest to the editors of the journal to reject the paper.}}'

The main judgement of the referee report `does not give any new insight' is not justified.
The report calls this {\it judgement} a {\it fact}, which is a logical fallacy.
In \cite[Remark 2.b]{AKS-1} we do explain why we do give a new insight:

`We think counterexamples of Theorem 1 are mostly interesting because their proof requires non-trivial ideas,
see below'.

%The main judgement of the report is {\it `does not give any new insight on \cite[Conjecture 5.5]{BFZ}'}.
%This is not justified; in the report a reference to a justification of the opposite
%(i.e., to \cite[Remark 2.b]{AKS-1}) is suppressed; see details in (b) below.

The referee suppresses the reference to this explanation in the above paragraph of the report, presumably because his judgement `does not give any new insight' does not withstand the test of being read next to this explanation.

\smallskip
{\bf (b)} The paragraph of the report quoted in (a) states that the rejection recommendation is based {\it only} on the (unjustified) judgement {\it `does not give any new insight'}.
After that, three pages of critical remarks are presented, starting with `About the evaluation of this paper'.
Thus the relation of these three pages to the rejection recommendation is unclear.

\smallskip
{\bf (c)} P. 1. `{\it Only recently, combining the work of Mabillard-Wagner \cite{MW14}
and Blagojevi{\'c}-Frick-Ziegler \cite{BFZ14}, Florian Frick in \cite{Fr15} presented counterexamples
to topological Tverberg conjecture in the case when n is an arbitrary non prime power.}'

%because the counterexamples were found in a series of papers by M.~\"Ozaydin, M. Gromov, P. Blagojevi\'c,
%F. Frick, G. Ziegler, I. Mabillard and U. Wagner.

This is misleading, and is publicly criticized in \cite[Remarks 1.9, 5.3 and \S5]{Sk16}
(most of the material of \cite{Sk16} cited here is available in earlier arXiv versions of this survey, and was even sent to the main `players of the game' before putting it to arXiv.)
For the reader's convenience, we quote \cite[Remarks 1.9]{Sk16} as Remark \ref{r:quot}.
Justified criticism of the description of references to the counterexample presented in \cite{Sk16} is not publicly available, see details in \cite[Remark 5.3]{Sk16}.

{\it Thus the report attempts to misuse the anonymous peer review system} to suppress the clear description of references concerning the topological Tverberg conjecture, presented in  \cite{Sk16}, presumably because that description is way too sound to be publicly criticized.
(See below description of further such attempts.)

%specific parts of \cite{Sk16} the referee disagrees with).
%In particular, in the above phrase Gromov's contribution is not mentioned,
%while Blagojevi{\'c}-Frick-Ziegler's [5] and Frick's [7] rediscovery  of the Gromov's implication
%`the topological Tverberg theorem, whenever available, implies the van Kampen-Flores theorem'
%is unduly named `presenting counterexamples combining other works'.

{\it We made the paper \cite{AKS-1} independent of that discussion}, by citing in \cite[Remark 2.a]{AKS-1} the surveys on the topological Tverberg conjecture containing different description of references, including \cite{BZ16} which contains a misleading description \cite[Remark 5.7]{Sk16}.
In the new version \cite{AKS-2} we added references \cite{BBZ, BS17, Sh18}, although they do not give anything new to that description comparatively to \cite{BZ16}, and contain misleading descriptions
\cite[Remarks 5.4, 5.5, 5.6]{Sk16}.

P. 5. '{\it the reference \cite{Sk16}: only published version should be quoted, otherwise this reference contains 4 or 5 different manuscripts,}'

It is usual to have many versions on arXiv and to cite a paper meaning the latest version.

P. 3. '{\it In Remark 2(a), which follows after Theorem 1 and offers motivation,
the authors, for the topic of the counterexamples to the topological
Tverberg conjecture, quote `surveys \cite{BZ16, Sk16}'. The proper way
would be to cite the original publications of the counterexamples \cite{Fr15}
(announcement) and \cite{BFZ} (journal version).}'

The papers \cite{Fr15, BFZ} are not the original publications of the counterexamples,
at least not the only ones, as explained in \cite[Remark 1.9, 5.4.abef, 5.5abcd, ]{Sk16}.

P. 4. '{\it In the paragraph after Remark 4, as an explanation for Theorem 5, counterexamples to the $r$-fold van Kampen-Flores conjecture are just mentioned and not explained.
Three papers are quoted in relation to this topic, ignoring the papers where the first counterexamples to the
$r$-fold van Kampen-Flores conjecture appeared.}'

The counterexamples are not just mentioned, but are explained, i.e., are stated together with proper references.
The referee does not name {\it `the papers where the first counterexamples to the $r$-fold van Kampen-Flores conjecture appeared'}, presumably  because if they would be named, it would be clear that it is unfair to name them as original references, see \cite[Remarks 1.9.b, 5.4.e, 5.5.b]{Sk16}.

We added in \cite{AKS-2} as a further explanation, the reference \cite[\S1, Motivation \& Future Work, 2nd paragraph]{MW14} to the extended abstract well-known to the people in the area and cited by the referee.

\smallskip
{\bf (d)} P. 2. `{\it This paper has only one main result stated in Theorem 1.}'

This statement is irrelevant to the rejection recommendation
(and is also wrong because there are \cite[Theorems 5 and 6]{AKS-1}).

\smallskip
{\bf (e)} P. 2. `{\it Theorem 6 is a particular case of \cite[Theorem 3.6 and paragraph afterwards]{Ba93}.}'

This is wrong.
Indeed, \cite[Theorem 3.6]{Ba93} takes a group $G$ from a certain class and proves that there exists \emph{some} representation $W$ of $G$, for which there exist $G$-equivariant maps $X \to S(W)$ for certain $G$-spaces $X$.
However, $G=\mathfrak S_n$ does not satisfy the conditions of \cite[Theorem 3.6]{Ba93}.
Also Theorem 6 states that for a specific group $G=\mathfrak S_n$ there exist $G$-equivariant maps from $X$ to a \emph{specific} representation of $G$, possibly not the one given by \cite[Theorem 3.6]{Ba93}.

\smallskip
{\bf (f)} P. 2. `{\it The result of Lemma 10 is known, ...}'

Since the report presents no reference to the result of \cite[Lemma 10]{AKS-1}, one has to treat this statement as incorrect.

P. 2. `{\it ...and can also be deduced from the existence of a $\Sigma_r$-map
$S^{2r-3}_{\Sigma_r}*\Sigma_r\to S^{2r-3}_{\Sigma_r}$...}'

%In fact, the referee does not present such an alternative proof:
The existence of a $\Sigma_r$-map $S^{2r-3}_{\Sigma_r}*\Sigma_r\to S^{2r-3}_{\Sigma_r}$ is trivially equivalent to non-trivial \cite[Lemma 10]{AKS-1}.
The referee neither proves the existence nor gives a reference for it.
This should be compared to the referee's critical remarks on \cite[proof of Lemma 10]{AKS-1}, see (w1) below.

%A lemma having an unpublished alternative proof does not make a lemma known uninteresting.

\smallskip
{\bf (g)} P. 3. `{\it The paper is written in a very unprofessional
and unclear way, with incomplete proofs, and missing relevant references.
The state of the paper under review was to me particularly puzzling
after considering the paper, referenced by \cite{AK19}, authored by Avvakumov
and Karasev (2/3 of the present authors). The difference in quality is stunning.}'

These judgements are not justified.
Neither it is explicitly written that they are justified by critical remarks presented below in the report,
presumably because it is too clear that they are not justified.
Observe that \cite{AKS-1} presents a detailed (however imperfect) proof of \cite[Lemma 10]{AKS-1}, comparatively to the reference to a picture in an analogous situation in \cite{AK19}.
Thus this critical remark of the referee is inconsistent with his/her criticism from (w1) below.
%This suggests
%shows that the referee praises/criticizes papers without thoroughly reading them, and misuses the anonymous
%peer review system in order to punish one of the authors for his criticism \cite[\S5]{Sk16}, \cite{Sk17o}, %\cite[\S5]{Sk17-3} (because the referee failed to question this criticism in an open discussion).

%{\bf Other replies to specific remarks by the referee.}
%\footnote{Recall that the remaining remarks are more technical and can be omitted by a generic reader.}

\small

\smallskip
{\bf (b1)} P. 2. `{\it First, Theorem 1 does not solve any case of the \cite[Conjecture 5.5]{BFZ} --- as the authors state in the sentence which breaks from the first to the second page of the manuscript.}'

%Many papers published in Israel J. Math.

Even if a paper does not completely solve any case of motivating conjectures or problems, the paper could be an interesting or significant contribution to them.
The referee does not justify that \cite{AKS-1} does not contain such a contribution.
See (a,b) above.

\smallskip
{\bf (c1)} P. 2. `{\it Second, in the case when $\frac{d+3}r < \lceil \frac{d+2}{r+1}\rceil$ the result of Theorem 1 is an elementary fact, any affine map in general position yields and almost n-embedding even without any
condition on $n$. ... In some other situations the gap given in \cite[Lemma 5.2]{BFZ} is greater than the one given by Theorem 1.}'

The existence of trivial / known particular cases does not contribute to making a theorem uninteresting.
See the following phrase at \cite[bottom of p. 1]{AKS-1}:

`Theorem 1 provides even stronger counterexamples to the topological Tverberg conjecture:
{\bf for $d$ large compared to $r$}...'

See also the phrase of \cite[Remark 2b]{AKS-1} cited in (a) above.

(Note that \cite[Lemma 5.2]{BFZ} does not explicitly state any counterexamples, so the phrase `the gap given in
\cite[Lemma 5.2]{BFZ}' is meaningless.
No examples of `is greater' is presented either in \cite{BFZ} or in the referee report.)

\smallskip
{\bf (d1)} P. 2. `{\it Lemma 9 is a special case of \cite[Lemma 3.9]{Ba93} where $Y$ is taken to be the empty set.}

P. 5. '{\it After Lemma 9, a reference towards \cite[\S5]{AK19} is given for the proof
and relevant information. As explained Lemma 9 is a special case of \cite[Lemma 3.9]{Ba93}.}'

We are grateful to the referee for bringing \cite[Lemma 3.9]{Ba93} to our attention.
We have added the following sentence to \cite{AKS-2}:

`In particular, this lemma follows from \cite[Lemma 3.9]{Ba93}, although
to read the direct proof in \cite[\S5]{AK19} is simpler than to find the notation required for
\cite[Lemma 3.9]{Ba93} and make such a deduction.'

Since \cite[Lemma 9]{AKS-1} is not a result of \cite{AKS-1}, but is attributed to \cite{AK19},
the above remark of the referee cannot contribute to the rejection recommendation.

\smallskip
{\bf (e1)} P. 3. `{\it I would like to point out that authors systematically quote their own papers
for the results of others, which cannot be seen as professional.}'

This is wrong, and is not justified by any reference to the paper under review.
%Neither it is written that the critical remarks presented below in the report justify this judgement,
%presumably because it is too clear that they do not.

For the counterexamples to the topological Tverberg conjecture, we quote the surveys  \cite{Sk16, BZ16}.
The surveys contain both original references and different description of them.
This is done in order to make this paper independent of the discussion of such a description, cf. (c) above.
We do not exclude from our citations the survey \cite{BZ16} written by other authors, and containing a misleading \cite[Remark 5.7]{Sk16} description.
% (attributing to \cite{BFZ, Fr15} the counterexamples to the topological Tverberg conjecture).
This cannot be called {\it `authors systematically quote their own papers for the results of others'}.

%!!! sent

\smallskip
{\bf (f1)} P. 3. `{\it The title of the paper says `Stronger counterexamples. . . '. It is unclear what do you mean.
As already explained one way to talk about this is to analyze the gap (1), but only as a function of `dimension'
and $n$.
Otherwise there is are new `stronger' counterexamples.}

What we mean is the usual `stronger than previously known'.
This is explained more precisely in the following sentence of \cite[Remark 2a]{AKS-1}:

`Theorem 1 provides even stronger counterexamples to the topological Tverberg conjecture:
for $d$ large compared to $r$ we have $N>(d+1)(r-1)$ and even $N>F$.'

The referee's phrase {\it ``Otherwise there is are new `stronger' counterexamples''} does not make sense to us.

\smallskip
{\bf (g1)} P. 3. '{\it Thanks should be given in a standard way, not by appending `*' to the title.}'

This could be fixed during the copy editing process if the paper is accepted.

\smallskip
{\bf (h1)} P. 3. '{\it The abstract should concisely present what is done in the paper and what is actually new...}'

This is done, and the referee does not justify that this is not done.
The Theorem stated in the abstract is new.
In \cite{AKS-2} we added the following phrase to the abstract:

`This was improved in 2015 by P. Blagojevi{\'c}, F. Frick, and G. Ziegler using a simple construction of higher-dimensional counterexamples by taking $k$-fold join power of lower-dimensional ones.'

\smallskip
{\bf (i1)} P. 3. '{\it ... In particular, please indicate precisely in what sense your counterexamples are new and stronger than done in previous work.}'

Since the relation involves technical formulas, we present it in \cite[Remark 2a]{AKS-1}, not in the abstract.

\smallskip
{\bf (j1)} P. 3. '{\it Please unify the notation used in the first paragraph of the paper,
the abstract and the statement of Theorem 1.}'

The referee indicates no conflict of notation, and we could not find it.

\smallskip
{\bf (k1)} P. 3. '{\it It is unclear what $K$ and $\Delta_N$ are.}'

This is explained in the first two sentences of the paper (see also Remark 4).

\smallskip
{\bf (l1)} P. 3-4. '{\it This presentation of fairly standard notions is very confusing.
In Matou\v sek's book \cite{Ma03} these notions are presented well.}'

The judgement of the first sentence is not justified.
The fact that some exposition is good does not show that an alternative exposition is `very confusing'.
We wrote the paper assuming that the reader of Israel J. Math. need not have studied  the book \cite{Ma03}.

% by any references to the paper under review.
%We did not find any confusing presentation of standard notions.
%, see above.
%Our exposition is a bit shorter, and so hopefully is more accessible to non-specialists.

\smallskip
{\bf (m1)} P. 4. '{\it In Remark 2(a), you talk about \cite[Conjecture 5.5]{BFZ}, but do not fully state it or explain it.}'

We do fully state the unknown part of \cite[Conjecture 5.5]{BFZ} in Remark 2a.
We omit the known part of the conjecture (which is thus not a conjecture, but a result),
i.e., the first line of \cite[Conjecture 5.5]{BFZ}.

\smallskip
{\bf (n1)} P. 4. '{\it In the first sentence of the second paragraph of Remark 2(a), you claim:
`This problem was considered in \cite[\S5]{BFZ} cleverly using the pigeonhole principle.'
This is inaccurate, the idea is different.}'

%The pigeonhole principle is used in [BFZ, \S5], in order to prove that the `$k$-fold joins' constructed there
%indeed provide stronger counterexamples.

In \cite[Remark 2a]{AKS-2} we changed our phrase to

`This problem was considered in \cite[\S5]{BFZ}, where higher-dimensional counterexamples were
constructed from lower-dimensional ones'

and added a short construction accessible to non-specialists.

\smallskip
{\bf (o1)} P. 4. '{\it In Remark 2(b), you explain why you do not compare the results you
prove with the already known results.}

This is wrong.
We do compare in \cite[Remark 2a]{AKS-1} the results we prove with the already known results.
In \cite[Remark 2b]{AKS-1} we explain a different thing:
`...We do not spell out even stronger counterexamples which presumably could be obtained  by combining the procedure of \cite[\S5]{BFZ} with Theorem 1.'

\smallskip
{\bf (p1)} P. 4. '{\it The argument says `We think counterexamples of Theorem 1 are mostly interesting because their proof requires non-trivial ideas, see below.'
I dislike the term `nontrivial idea' which is subjective, and implicitly implies that previous work is trivial.}'

The words `nontrivial idea' are common usage in introduction to research papers.
For a reader, it is important to know the authors' understanding of what is simpler and what is more complicated.
It is correct that the proof of Theorem 1 of \cite{AKS-1} is more complicated than \cite[Lemma 5.2]{BFZ},
see a short exposition of \cite[Lemma 5.2]{BFZ} accessible to non-specialists in \cite[Remark 2a]{AKS-2}.

The referee did not explicitly state the opposite (only hinted that), presumably because the opposite does not withstand the test of being explicitly stated and looked upon.

\smallskip
{\bf (q1)} P. 4. '{\it In Remark 4 (motivation) a definition of a simplicial complex and its geometric realization is given. First, a definition is called `Remark' and in addition `(motivation)' is added for complete confusion.}'

Remark 4 of \cite{AKS-1} is not a definition but is a remark containing definitions, a convention, an observation, and information.
%So in order to avoid confusion we called it `Remark' not `Definition'.
The word `(motivation)' is proper because this remark motivates \cite[Theorem 5]{AKS-1}.

\smallskip
{\bf (r1)} P. 4. '{\it Next, a standard term of `geometric realization' is substituted with the term `body',
and then nowhere in the paper this term reappears.}'

We mention both terms: `The {\it body} (or geometric realization)...'.
This is convenient to readers working in different areas who might be more familiar with one or with the other term.
This is unimportant because indeed nowhere in the paper this term (or this notation) reappears,
which is explained in the last phrase of the remark:

`We abbreviate $|K|$ to $K$; no confusion should arise.'

\smallskip
{\bf (s1)} P. 4. '{\it Reference to \cite[Theorem 3.3]{Sk16} is inappropriate since this is just the original statement of \"Ozaydin result.}'

This is wrong because the paper \cite{Oz} does not contain the statement \cite[Theorem 3.3]{Sk16}
(although there is no question of attributing it to \"Ozaydin, which is done in \cite{Sk16}).
Reference to \cite{Sk16} is the more appropriate because this survey provides a simple and accessible to non-specialists proof of the \"Ozaydin result.

\smallskip
{\bf (t1)} P. 4. '{\it ... but with incorrect details in the sketch of the presented proof \cite[\S3.2]{Sk16}.}'

The referee's statement that there are incorrect details in a proof (in a published paper) is inappropriate
because the referee does not indicate which phrase(s) of \cite{Sk16} exactly is (are) incorrect.

\smallskip
{\bf (u1)} P. 4. '{\it There is no reason to mention \cite{BG17} at all.}'

We do give the reason in the phrase mentioning \cite{BG17}.
The referee does not explain what is wrong with it.
The suggestion to remove the mention of \cite{BG17} is inconsistent with the referee being touchy about not citing other papers, and is a minor issue which cannot contribute to the rejection recommendation
(we agree to remove the mention of \cite{BG17} if the referee or editors insist).

\smallskip
{\bf (v1)} P. 4. '{\it The statement of Theorem 7 quotes three papers without proper credit to precedence.}'

The statement of \cite[Theorem 7]{AKS-1} quotes three papers in alphabetical order.
The referee does not explain what is wrong with that, and what he/she means by `precedence',
presumably because his/her objection does not withstand the test of being explicitly stated and looked upon.

The following historical details are not used in \cite{AKS-1}, and so are not presented there.
The paper \cite{Sk17o} (and \cite[\S5]{Sk17-3})

$\bullet$ contains a Skopenkov-Wagner description of the problems with the proof \cite{MW16};

$\bullet$ attributes the result claimed in \cite{MW16} to Mabillard-Wagner, in order to concentrate on mathematics, not on priority questions.

\smallskip
{\bf (w1)} P. 5. '{\it Proof of Lemma 10 is very hard to follow, especially in the final part `construction of $f'_+$'.
You are talking about homotopy between a map with degree -1 and degree 1 --- more care is needed.
Then `the equivariant Borsuk Homotopy Extension Theorem' is used, with no
reference, no assumptions checked for its application.
Thus, the proof of Lemma 10 is on the level of a sketch.}'

The second sentence is wrong, we do not talk about homotopy between maps of degree $-1$ and degree 1.
The referee does not refer to specific sentence where we do (because there are none).
(E.g. before the second display formula in p. 4 we write that the map between spheres is defined for $t=0,1$, not for $t\in[0,1]$.)

The third sentence is a valid criticism.
However, the referee indicated no significant problem, and this criticism is easy to take into account
\cite[pp. 5-6, Proof of the assertion: construction of $f_{1,+}$]{AKS-3}.

Thus this paragraph of the report cannot contribute to the rejection recommendation.
%Cf. (b) above.

%The referee does not justify that the proof is very hard to follow because of unclear exposition,

%We do provide necessary care by writing that the range of the homotopy $h$ is Eucleadean space,
%while the range of maps $f,f_0,f_1,f'_+,f'_-$ (whose degree is mentioned) is sphere.

%omit??? In our opinion, the theorem is a textbook result and so could be used without reference.

%The referee's justification of this conclusion is wrong, see just above.
%Thus conclusion should be compared with (g) above.

\smallskip
{\bf (x1)} P. 5. '{\it For the references:

- the reference \cite{BFZ} should be updated, the paper is published,

- the reference \cite{BG17} should not be used,

- the reference \cite{BZ16}, the last name of the first author is misspelled, the year of publication and the name of the publisher are missing,

- the reference \cite{Fr15} should be updated, number, year, and pages are missing,

- the reference \cite{Oz}, the year is missing,

- for the references \cite{AK19, MW15, MW16, Sk17o} the years of posting should be added as well as the number of pages, even some date can be deduced from the arXiv reference of a choice of quotation label.}

This could be fixed during the copy editing process if the paper is accepted.
However, we did update the references in \cite{AKS-2}.

\smallskip
{\bf (y1)} P. 5. '{\it At the end, in the last sentence of the paper: `*' should not be used to classify the papers in the situation when all the relevant surveys and/or books are not quoted.

Some of the survey papers and the books related to the topic which
are not quoted are \cite{Ma03, Zi11, BBZ, BS17, Sh18}.}'

We do cite the surveys on the topological Tverberg conjecture containing different description of references, including \cite{BZ16} which contains a misleading \cite[Remark 5.7]{Sk16} description.
We did not cite the surveys \cite{BBZ, BS17, Sh18} because they do not present any new description comparatively to \cite{BZ16}.
However, we added citations of the surveys \cite{BBZ, BS17, Sh18} in \cite{AKS-2}.

We did not cite (and did not add to \cite{AKS-2} citations of) the book and the survey  \cite{Ma03, Zi11} because they do not mention counterexamples to the topological Tverberg conjecture (these references appeared before the counterexamples).
So these nicely-written references are too far afield to be mentioned in a research paper on stronger counterexamples.
\normalsize
\end{example}

\begin{remark}\label{r:quot} This is a quotation of \cite[Remarks 1.9.abc]{Sk16}.
We omit the footnotes containing details of justification; the numbers refer to numbered statements of \cite{Sk16}.

\smallskip
\textbf{(a)} {\it `The topological Tverberg theorem, whenever available, implies the van Kampen-Flores theorem'} \cite[2.9.c, p. 445, lines --1 and --2]{Gr10}.
This is the Constraint Lemma 1.8.
What M. Gromov called `theorem', we call `conjecture'.
See (c) for the rediscovery of this lemma.

This lemma was proved in \cite[2.9.c, p.446, 2nd paragraph]{Gr10} by a beautiful combinatorial trick reproduced in the current paper.

\smallskip
\textbf{(b)} {\it A counterexample to the $r$-fold van Kampen-Flores conjecture for $r$ not a prime power}
(Theorem 1.7) was implicitly obtained by I. Mabillard and U. Wagner in \cite{MW14, MW15} by mentioning the \"Ozaydin  Theorem 3.3, and proving that the \"Ozaydin  Theorem 3.3 implies Theorem 1.7.
Indeed,

$\bullet$ the \"Ozaydin  Theorem 3.3 states that some equivariant map exists for $r$ not a prime power;

$\bullet$ \cite[Theorem 3]{MW14} states that for codimension at least 3 the existence of the above equivariant map implies the existence of an almost $r$-embedding as in Theorem 1.7.

Failure of the $r$-fold van Kampen-Flores conjecture for $r$ not a prime power was first {\it explicitly stated} by F. Frick \cite{Fr15}.

\smallskip
\textbf{(c)} {\it Since $O$ and $O\Rightarrow\overline{VKF}$ and $TTC\Rightarrow VKF$ imply $\overline{TTC}$,
by (a) and (b) papers \cite{Oz, Gr10, MW14} together give counterexample to the topological Tverberg conjecture} (so that there remained the non-trivial task of writing versions of \cite{Oz, MW14} that could be accepted
to
%ABS: by?
peer-review journals).
However, neither Mabillard and Wagner nor the topological combinatorics community were aware of (a) before 2016.
This is surprising because the next part \cite[2.9.e]{Gr10} of Gromov's paper was discussed during the problem session at 2012 Oberwolfach Workshop on Triangulations.
So this community saw a serious problem with the Mabillard-Wagner approach to a counterexample to the topological Tverberg conjecture: maps from the $(d+1)(r-1)$-simplex to $\R^d$ do not satisfy the codimension $\ge3$ restriction required for \cite[Theorem 3]{MW14}; these maps actually have negative codimension.
F. Frick \cite[proof of Theorem 4]{Fr15} realized that this problem can be overcome by (a).
He did this by rediscovering (a), not by finding the reference \cite[2.9.c]{Gr10}.
In fact, (a) was implicitly rediscovered earlier by Blagojevi\'{c}-Frick-Ziegler
\cite[Lemma 4.1.iii and 4.2]{BFZ14}.
\end{remark}

\begin{example}[the report on handling of \cite{AKS-2} in Proc. of the Amer. Math. Soc.]
%  20.03.2022
\label{r:akspamsd} (R. Karasev and A. Skopenkov)
The paper was rejected on the basis of a referee report X shown to be incompetent below.

\smallskip
\textbf{(a)} The referee's specific critical remarks are either

$\bullet$ attempts to misuse the anonymous peer review system to promote an opinion which failed to be substantiated by an open discussion, see (b) below, or

$\bullet$ wrong or unjustified statements and judgements, see (c)-(f) and (h) below, or

$\bullet$ valid critical remarks that cannot contribute to the rejection recommendation, see (g) below.

% Namely, to suppress clear description of references concerning the topological Tverberg conjecture [Sk16]  %(instead of publicly criticizing specific parts of [Sk16] the referee disagrees with).
%Together with the incompetence of other remarks, R(1) that the report misuses the anonymous peer review system
%to promote an opinion which does not withstand an open discussion.

The most important critical remarks to the report are (d), (e), and (h) below.
%The referee remarks are given below in italics and in quotation marks.
%Both the referee report and this reply refer to numbered statements and references from \cite{AKS-2}.
%arXiv:1908.08731v2.

\smallskip
\small
{\bf (b)} {\it `... Recently, in 2015, Frick presented counterexamples to the topological Tverberg conjecture
for non prime powers...'.}
(We do not quote from the report the rest of the introductory paragraph.)

In  Example \ref{r:aks-ijm}.c
%5.7.c of arXiv:2101.03745
it is explained that {\it this is misleading}, and that {\it the paper under review is independent of
the discussion of references on the counterexample to the topological Tverberg conjecture}.

\smallskip
{\bf (c)} {\it `The paper under review studies further the `dimension gap' of counterexamples.
The main and the only result of the paper is \cite[Theorem 1]{AKS-2} which gives
improved ``dimension gap'' of counterexamples to the topological Tverberg conjecture.'}

Here {\it `the only'} is incorrect (because the paper has \cite[Theorems 4 and 6]{AKS-2})
and is irrelevant to the rejection recommendation.
Cf. Example \ref{r:aks-ijm}.c.
%5.7.c of arXiv:2101.03745.
\normalsize

\smallskip
{\bf (d)} {\it `The authors do not properly evaluate and compare their result
to the other known related results. The claim in \cite[Remark 2.a]{AKS-2} ``Theorem 1
is a partial result on \cite[Conjecture 5.5]{BFZ} stating that ...'' is not correct.
Indeed, if \cite[Theorem 1]{AKS-2} would address any instance of \cite[Conjecture 5.5]{BFZ} the result
would definitely deserve publications in a good journal.'}

Even if a paper does not completely solve any case of motivating conjectures or problems, the
paper could be an interesting or significant contribution to them.
The referee does not justify that the paper does not contain such a contribution.
Cf. Example \ref{r:aks-ijm}.b1
%5.7.b1 of arXiv:2101.03745,
as well as the above quotations from report Z and quick opinion Y.
(E.g. `$f(n)\ge n-1$' {\it is a partial result} on a conjecture `$f(n)=n^2$', even though
$f(n)\ge n-1$ does not imply that $f(n)=n^2$ for any instance of $n$.)
Thus the claim in \cite[Remark 2.a]{AKS-2} is correct, and the referee's conclusion {\it `do not properly evaluate and compare...'} is incorrect.
Still, we do not object to slightly changing the correct phrase of \cite[Remark 2.a]{AKS-2} to emphasize that our
partial result on \cite[Conjecture 5.5]{BFZ} does not prove the conjecture for any instance.

\smallskip
{\bf  (e)} {\it `The rest of the paper gives the proof of \cite[Theorem 1]{AKS-2} is a way which is very hard to follow and even
harder to verify.'}

This conclusion is not justified by references to the paper and suggestions what could be done in a more clear way, except for a single minor issue (see below our replies to the following phrases).
So the report does not show that the paper is hard to read because of poor exposition, not just because of complexity of the matter.
Cf. the above quotations from report Z.
We wrote this paper in a clear way according to the high standards of \S\ref{s:differ},
%arXiv:2101.03745,
including sending it to colleagues and working on their critical remarks before journal (and even arXiv) submission.

\smallskip
\small
{\bf  (f)} {\it `The main technical ingredient of the proof of \cite[Theorem 1]{AKS-2} is  \cite[Lemma 10]{AKS-2}. It constructs degree zero equivariant self-maps of odd dimensional spheres.
The proof of \cite[Lemma 10]{AKS-2} is not well presented.
It is given in a paragraph on \cite[page 4]{AKS-2}, but relies on a sequence of claims which are presented on the remaining pages of the paper.'}

The last sentence describes what we believe a {\it well-structured proof} is.

\smallskip
{\bf  (g)} {\it `A problematic issue occurs is the use of the so called ``equivariant Borsuk Homotopy Extension Theorem''.
The authors give no reference to such a theorem but instead give a comment in \cite[footnote 3]{AKS-2} by stating the needed result and claiming ``The proof is analogous to the non-equivariant version \cite[\S5.5]{FF89}.'' The statement given in \cite[footnote 3]{AKS-2} is obviously not true.'}

We are grateful to the referee for finding a flaw in our argument.
This flaw is easy to fix, see the next version \cite{AKS-3}.
%arXiv:1908.08731v3.
Of course, a referee need not think how to fix a flaw he/she found.
However, a competent refereeing includes revision not rejection recommendation in the case of a flaw in the proof not shown to be serious by the referee.
\normalsize

\smallskip
{\bf  (h)} {\it ` In summary, the paper considers an interesting and relevant problem in Topological Combinatorics.
It presents a partial results on an improved ``dimension gap'' of counterexamples to the topological Tverberg conjecture, but fails to give a comprehensive and complete proof of the result.
For these reasons I suggest to the editor of the PAMS to reject the submitted version of the
paper ``Stronger counterexamples to the topological Tverberg conjecture''.
The only suggestion I can make to the authors is to rewrite the paper completely
with much more care for the reader, and, more important, for the arguments in the proofs.'}

Here the judgement {\it `fails to give a comprehensive and complete proof of the result'} and the suggestion
{\it `to rewrite the paper completely'} are unjustified and unsuitable, as justified by (c)-(g) above.
% (and the revised version arXiv:1908.08731v3).
\end{example}

\begin{example}[the report on handling of \cite{KS21e} in Fund. Math.]\label{r:ksfm}
The paper \cite{KS21e} (shortened by deleting sections 2.4.and 2.5) was submitted in April 2022.
It was rejected in December 2023 because of

$\bullet$ a logical fallacy in the Editor's letter:
from the presence of earlier published results (or combinations of those) in the introduction, the letter concludes (ignoring the new results) that `\emph{we are close to ... the main results having already been obtained in other publications}'.

$\bullet$ a wrong statement in the Editor's letter: `\emph{some of the results of the removed Appendix of \cite{KS21e} still play a role in the proofs of the submitted paper}'.

See details in (a) below.
No other criticism of the paper was presented by the Editor.
No referee reports were presented by the Editor.
%See also Remark \ref{r:ksfm-l}.

\smallskip
\small
\textbf{(a)}
The short version of the reports presented in the Editor's [in the letter to the Editor: in your] letter

* states `the Authors themselves say rightly that several of their main results are a combination of those which appeared already in [PT19], [Me06], [FKT], or [AMS+]';

* ignores the novel part, and

* concludes that `we are close to' `the main results having already been obtained in other publications'.

I'm afraid this is a logical fallacy.
Presence of earlier published results (or combinations of those) in the introduction does not mean that the paper does not have new results.
(E.g. the fact that Theorem 1.3.1 is attributed to [vK32, Sh57, Wu58, Me06, FKT, AMS+] cannot contribute to the lack of novelty in other results.)

The paper does not say that several of its main results are a combination of those which appeared earlier.
Results which are combinations of those which appeared earlier are presented as earlier known results, not as our results: we give clear references to earlier publications.
(See e.g. reference to Remark 1.3.7 in p. 1, and in Comments on the proof after Theorem 1.3.1: ` Part (a) is essentially classical ...  Part (b) is proved in ...  Part (c) is proved in ...'.)

The short version of the reports presented in your letter does not indicate any specific place where results of the removed Appendix of \cite{KS21e} are used in the proofs of the submitted paper.
In my opinion, there are no such places; the appendix is only mentioned for introductory purposes.
\normalsize
\end{example}

\begin{example}[the report on handling of \cite{DS-3} in Intern. Math. Res. Notes]\label{r:dsirmn}
The paper \cite{DS-3} (shortened by deleting \S3, \S4 and footnotes) was submitted in February 2024.
It was rejected in March 2024 based on an incompetent report, as explained in (a,b,c) below.\footnote{\label{f:dsirmn} In reply to the letter of Remark \ref{r:dsirmn-l}.b and Example \ref{r:dsirmn}.abc, the Editor suggested submission of the new version, without withdrawing the rejection decision.
However imperfect, this was a reasonable compromise.
So Example \ref{r:dsirmn} does not contribute to the judgment on Intern. Math. Res. Notes presented here.
However, Example \ref{r:dsirmn5} does.
Besides, later development described in Example \ref{r:dsirmn5}.c9 shows that the referee (although not the Editors) of \cite{DS-3} did violate scientific ethics, by misusing the anonymous peer review system to suppress the publication of a paper in favor of another competitive paper that appeared later.}

\smallskip
\textbf{(a)} The report is incompetent, since

$\bullet$ it states that the main results of the paper under review is known without giving any specific reference;

$\bullet$ the statement that the PL case of the problem (i.e. of the K\"uhnel conjecture on embeddings)
is solved using hard Lefschetz property disagrees with existing references (see explanation [added later: in (c)] below).

\small
\smallskip
\textbf{(b)} Less importantly, the report ignores

$\bullet$ the equivalence of PL and topological versions of the K\"uhnel conjecture for $k > 2$, which is well-known, see the Patak-Tancer paper [arXiv:1904.02404, Proposition 7] cited by the referee
(see also \cite[the second paragraph of Remark 1.3.c]{DS-5});

$\bullet$ the version for $\Z_2$-embeddings proved in Theorem 1.8 of the paper under review, which version implies the topological version (the implication is simple and known, see \cite[Remarks 1.3.c and 1.9.a]{DS-5}).

\smallskip
\textbf{(c)} The only paper mentioning the K\"uhnel conjecture on embeddings in connection with hard Lefschetz property is the unpublished paper \cite{Ad18}.
In \cite[Section~1.6, (1)]{Ad18} one concerns a different K\"uhnel conjecture [K\"uhnel 94, Conjecture C] [added later: this is reference [Ku94] of \cite{DS-3} or of \cite{DS-5}].
In \cite[Remark~4.9]{Ad18} one reads:
'if a complete $k$-dimensional complex on $n$ vertices embeds into $M$ sufficiently tamely (so that it extends to a triangulation of $M$), then $\binom{n-k-1}{k+1} \le \binom{2k+1}{k+1} b_k(M)$'.
This is the K\"uhnel conjecture on embeddings under additional assumption
`the embedding extends to a triangulation of $M$'.
It is unknown whether this additional assumption can always be achieved (see details in \cite[Remark 1.3.g]{DS-5}).
Thus the K\"uhnel conjecture on PL embeddings remains unsettled.

Less importantly, the proof under the additional assumption in \cite[Remark~4.9]{Ad18} uses another conjecture \cite[Remark 4.1.a]{DS-5}.
Criticism of the claim for this conjecture is publicly available \cite[footnote 13]{DS-5}.

For more details on references to the K\"uhnel conjecture on embeddings see \cite[Remark 1.3]{DS-5}.
\normalsize
\end{example}

%(b) The report. [Available upon request]
%There are two versions of the problem: the one considered by Patak and Tancer (they work in the larger topological category), which is more general and still open.
%And the one considered in this paper (provided the maps are PL), which is solved.
%The ingredients are as follows: Hard Lefschetz property implies Kühnel conjecture was proved by Gil Kalai, and later with a simpler proof by Adiprasito.
%That the hard Lefschetz property holds in the relevant case also has several proofs, by Adiprasito (somewhat long and complicated proof) and then two simpler proofs (Papadakis-Petrotou and Karu-Xiao).

\begin{example}[the report on handling of \cite{DS-5} in Intern. Math. Res. Notes]\label{r:dsirmn5}
The paper \cite{DS-5} (shortened by deleting \S3, \S4 and footnotes) was submitted in May 2024.
It was rejected in September 2024 based on an incompetent report, as explained in (a) below.

\smallskip
\textbf{(a)} The report by Referee \#1 of
%[added later: the paper under review]
[DS5] is incompetent, since

$\bullet$ it is partly based on factual mistakes (see (c5,c8) below);

$\bullet$ it leaves its crucial judgments without justification, and ignores the information in \cite{DS-5}
justifying the opposite judgments, instead of objecting to the information (see (c1) below);

$\bullet$ it criticises the exposition of \cite{DS-5} for its advantages, e.g. for its being accessible
to mathematicians not specialized in the area (see (c4,c11) below);\footnote{Added later: thus the two reports contradict to each other, cf. (c4, c11) to (d1-d3).}
%???

$\bullet$ in violation of scientific ethics, it misuses the anonymous peer review system
to suppress the publication of a paper in favor of another competitive paper
that appeared later (see (c9) below).

The report by Referee \#2 (not only is positive but also) cannot contribute to the rejection decision.
We are grateful to Referee \#2 for suggesting modifications.
We incorporated all but two minor comments (and explained in (d5,d6) what is wrong with the remaining two).
However, modifications are easy to do, not major, as we justify in (d1--d6).

See also other critical remarks on the report  \#1 in (c1--c12).
Most of the justification in (c1--c12,d1--d6) is accessible to mathematicians not specialized in the area.

\smallskip
\textbf{(b)} In (c1--c12) and (d1--d6) both the reports and our replies refer to the numeration in \cite{DS-5} (not in the latest arXiv version), which is the same as the numeration in the resubmitted version [DS5]
(if the numbered item is present in the resubmitted version)
[added later: (the numbered items discussed below are present in the resubmitted version)].
%At the beginning of each item a quotation from referee report is given in italics.

%%%%%%%%%%%%%%%%%%%%%%%% \#1

\smallskip
\textbf{(c1)} \emph{`The results are correct, but incremental and the proof method is just an extension of the previous proof.'}

The judgment `incremental' is unjustified because it ignores the fact that for 6 years
after the appearance of a linear estimate, no estimates stronger than linear ones appeared.

This report does attempt to justify its other judgments.
Those justifications are shown to be incorrect.
So those other judgments are poorly thought over or biased.
Together with the above, this shows that the judgment `incremental',
which the report even does not attempt to justify, is poorly thought over or biased.

Without analysis of a new proof, its being an `extension of a previous proof'
does not imply its triviality.
Almost any proof is `just an extension' of previous arguments,
including proofs of Perelman and Wiles.
Such an analysis is not provided in the report.
Moreover, \cite[\S1, subsection `Topological and linear algebraic parts']{DS-5}
does explain the non-triviality of the proof:
`Our theoretical achievement allowing us to prove Theorems 1.1 and 1.2
is to fit what we can prove in topology to what is sufficient for algebra.
Thus our main idea is the notion of an $(n,k)$-matrix.
Thus the proof is split into two independent parts.'
The report does not object to this description of non-triviality, but ignores it.

The report's remarks (given below) on the exposition in \cite{DS-5}, together with our comments
on these remarks, show that in fact the exposition is reasonably clear and well-structured.
Observe that usually a proof seems to be simple because of
a clear and well-structured exposition, not because of triviality and lack of novelty.

In addition, the judgment `incremental' disagrees with the judgment from report \#2:
`...the authors’ progress is a welcome advance on the status of this problem ...
the results are new and on the level to be published at IMRN'.

\smallskip
\textbf{(c2)} \emph{`The paper is very poorly written and most readers would be discouraged by its style.
... The main problem is that the paper is very unstructured,...'}

This is unjustified and wrong, as we explain below following the justification by the referee.

\smallskip
\textbf{(c3)} \emph{`...remarks are used throughout the text and serve various purposes: some of them are historical and serve as a motivation (e.g. Remark 1.3), some others are used to state the main result(!) (see Remark 1.3c), some other are used to state conjecture (Remark 1.3f).'}

There is nothing wrong (and the referee does not explain what is wrong) with remarks being used throughout the text and serving various purposes different from stating the main result.
The report wrongly states that Remark 1.3c is the main result, see (c4).

\smallskip
\textbf{(c4)} \emph{`The main theorem is Theorem 1.2, which immediately yields Theorem 1.1.
Remark 1.9 and Theorem 1.8 then show Theorem 1.2 holds in much greater generality for $\Z_2$-embeddings.
Why don't the authors want to state the result in full generality in the introduction?
In the abstract, the only claim the result for embeddings of skeleta of simplices to connected sums of $S^k \times S^k$?
It would be much better to state the full result and only mentioned Kühnel's conjecture as a consequence.'}

The result in full generality is stated in the introduction,
see subsection `Almost embeddability and $\Z_2$-embeddability'.
The name of the journal is `International Mathematics Research Notices',
not `International Research Notices on Embeddings'.
So we must make the introduction accessible to non-specialists in the area
(who are not familiar with almost embeddings and $\Z_2$-embeddings).
In particular, an earlier statement of the result in full generality would make the text
hardly accessible to referee \#2, who asked to define even as basic notions as
`intersection form' and `general position'.

\smallskip
\textbf{(c5)} \emph{`The crucial definition of an (n,k)-matrix is not stated as a definition,...'}

This is wrong, see the sentence after Lemma 1.7:
`A symmetric, independent, additive, non-trivial matrix is called an (n, k)-matrix.'

\smallskip
\textbf{(c6)} \emph{`...but the reader has to put it together from several pieces spread out on page 6.'}

There is nothing wrong (and the referee does not explain what is wrong)
with some definitions being used in another definition.

\smallskip
\textbf{(c7)} \emph{`The statements of other claims are unclear, without using external references,
e.g. P8L14: We did not suceed in generalizing to higher dimensions [DS22, Lemma 3.4] ...
Here the authors claim that they did not manage to generalized the lemma to higher dimensions;
but they do not bother to write down the precise statement ...
Moreover, this is in a part that should explain the proof idea.
I think this clearly illustrates, why the paper is poorly written.'}

\emph{`P8L14: We did not succeed in generalizing to higher dimension [ref], which implicitly appeared in []...
%->
$\to$
why are the readers forced to search external references?'}

We did bother to write down the precise statement of \cite[Lemma 3.4]{DS-5}, although in the arXiv version of the paper, not in the submitted version.
So it is easy to include the statement in the full version as required by the report
(and we did it; we give other external references only in remarks which,
as we explain, could be ignored).

This and the remarks to P2L19, P2L30--32, P6L24 below are the only proper objections
by the report \#1 (as justified by comments to other objections), and are minor ones.
So they cannot contribute to the rejection decision.
We think this clearly illustrates why the report is poorly written.

\smallskip
\textbf{(c8)} \emph{`Moreover, assumptions of the statements and claims have to be guessed on many places.
For example Conjecture 1.3.f. Is this conjecture stated for PL-manifolds or topological manifolds in general?'}

In \cite[the list of `Notation and conventions' before Theorem 1.1]{DS-5} one reads
'We consider only piecewise linear (PL) $2k$-manifolds.
Unless otherwise specified, we consider only PL maps'.

\smallskip
\textbf{(c9)} \emph{`For PL, it [Conjecture 1.3.f] has been proven in https://arxiv.org/abs/2404.12265, for non-triangulable manifolds it makes no sense.'}

Just as the paper
[AP24] := arXiv:2404.12265v1 (or v2, v3),
the report hides the fact that the paper [AP24] implicitly claims a competitive result.
This is explained in \cite[Remark 1.3.g]{DS-5}:
`The K\"uhnel conjecture on embeddings (a) follows from Conjectures [1.3.](d,f)'.
The paper arXiv:2404.12265 claims \cite[Conjecture 1.3.f]{DS-5}, while the paper \cite{Ad18}
%[Ad18] := arXiv:1812.10454v4
claims \cite[Conjecture 1.3.d]{DS-5}.
The latter claim is not withdrawn, although since April 2024 there is
an unanswered public criticism of that claim, see v4 or higher of \cite{DS-5} (footnote 16 in v4 and Remark 6.2 in \cite{DS-7}).

This implicit claim of the competitive result did not appear before the March 2024 report to \cite{DS-5} = ID IMRN-2024-227.
The paper [AP24] appeared on arXiv in April, 2024 after a post-referee-report discussion of our submission \cite{DS-5} involving one of the authors of [AP24], see arXiv:2208.04188v4, footnote 3.
The main results of \cite{DS-3, DS-5} are publicly available on arXiv since 1 Sep 2022
(v2 corrects a minor gap in the proof from v1 of 5 Aug 2022).
The paper \cite{DS-5} was submitted on 24 Feb 2024.
So scientific ethics require that no result could affect the acceptance/rejection decision
of \cite{DS-3, DS-5} if the result is claimed either by putting it on arXiv after 1 Sep 2022,
or by submission after 24 Feb 2024.
%In violation of that, the report \#1 brings the paper [AP24] into report to [DS5].
%So we would also be grateful if you could check if the referee \#1 of [DS5] has no conflict of interest in refereeing our paper.

(Perhaps bringing a later competitive paper to the report \#1 is not malevolence,
but the result of poor reading of \cite{DS-5}.
Such a difference can affect the Editors' moral decision on the immortal soul of the referee, but cannot affect the professional decision on the competence of the report \#1.)

\small

\smallskip
\textbf{(c10)} \emph{`P2L19: References to remarks could be ignored
%-->
$\to$ remarks could be ignored? It makes no sense to ignore just the references and not the remarks themselves.'}

\emph{`P5L30-32: The notation should be clarified explicitly. The part "let $a,a'$ be vertices from $K_{n,n}$ from different parts" makes no sense.
The authors probably mean if $a$ is a vertex in the first part, let $a'$ mean its counterpart in the other part (but then they need to specify what a counterpart is), ...'}

\emph{`P6L24: It is obvious that $A(f)$ is independent. --- This can be only said if $f$ is an embedding. (Otherwise, projective space provides a counterexample.) However, it is not clear, whether this is assumed at this point.'}

This is easily corrected; see the 2nd paragraph of (c7).

\smallskip
\textbf{(c11)} \emph{`P2L23: The analogues of our results are correct for topological embeddings (Remark 1.3.b), for almost embeddings (defined and discussed in Remark 1.3.c) and for $Z_2$ embeddings (defined and discussed in Remark 1.9). %-->
$\to$ This is an example how bad the paper is written. Instead of stating the main result in full form, the reader has to sieve through countless remarks to see what the result actually is.'}

This is an example of how reader-friendly the paper is written, not how badly the paper is written.
See justification in (c4) above.

\smallskip
\textbf{(c12)} \emph{`P2L53: see [KS21] for a simpler exposition
%-->
$\to$ the word simpler is questionable. In the original paper PT19, its authors gave a precise combinatorial description for embeddability of any k-complexes into PL $2k$-manifolds (assuming that either the complex or the manifold is sufficiently nice). Then, as an example of application they proved the linear bound for K\"uhnel's conjecture. In KS21, the bound is proved directly, but it is questionable whether it is easier to read.'}

In \cite[Remark 1.2e]{KS21},  it is explained why the direct proof of \cite{KS21} is easier to read (we added
[added later: reference to] \cite[Remark 1.2e]{KS21} to our phrase).
The report does not question that explanation.
Thus `easier to read' remains unquestioned.

%%%%%%%%%%%%%%%%%%%%%%%% \#2

\normalsize

\smallskip
\textbf{(d1)} \emph{`The current manuscript reads like an addendum to a series of papers, many unpublished.'}

The current manuscript uses not unpublished papers, but known results.
In some cases this was indeed not clear from our text, and we are grateful to the referee \#2 for observing this.
Such flaws are easily corrected, see below.

\smallskip
\textbf{(d2)} \emph{`For example, the proof of Theorem 1.8 is based on results from reference [IF].
This is not a peer-reviewed reference, it is a community-based wiki page which has
as a warning at the top of the page “this page has not been refereed.
The information given here might be incomplete or provisional”.
The authors should rewrite their proofs to make their manuscript not rely on such sources.'}

We just added (in \cite{DS-7})
%arXiv:2208.04188v7)
references to books and published papers.
Unless we receive an objection, in order to make the paper more self-consistent [added later: self-contained],
we plan to add section 5 of the arXiv v7 version into the resubmitted version
(although this is not directly required by the referee).

\smallskip
\textbf{(d3)} \emph{`Even for the definitions of general position the reader is sent to other unpublished papers.'}

We just added (in \cite{DS-7}) the definition, although the notion of general position is well-known, and so is used without definition in many papers in the area.

\smallskip
\textbf{(d4)} \emph{`If the authors intend for the remarks to be read at a second pass on the paper, they should move them to a section at the end of the paper, or it makes for cumbersome first reading.'}

We just added (in \cite{DS-7}) the following before Theorem 1.1:
`In this text remarks and references to them could be ignored by some readers,
but could be important for others; it is for a reader to choose whether to read a remark.'
We also just moved some important remarks to a separate section 2.

\smallskip
\textbf{(d5)} \emph{`Page 3, line 47 - This restatement cannot be correct.
Given a triangulation of said sphere, we can take barycentric subdivisions until the inequality stops being satisfied.'}

This argument suggested by the report \#2 does not work.
In the restatement the subcomplex is isomorphic (not just homeomorphic) to the complex. Barycentric subdivisions give only a homeomorphic subcomplex.

\small
\smallskip
\textbf{(d6)} \emph{`Page 4, line 10 - The Chicago manual of style indicates that when used in plural, `theorems' should not be capitalized.  Please homogenize throughout the paper.'}

We would be grateful if referee \#2 could indicate the specific point
from the Chicago manual of style concerning this situation.
For us `Theorem 3' but `theorems 3 and 5' in the same text look inconsistent.
We haven't seen any paper by a native English speaker styled like
`Theorem 3' but `theorems 3 and 5'.
(BTW, GPT chat suggests `The Chicago Manual of Style'.)
In any case these things could be done later if the paper is accepted.
\normalsize
\end{example}

\begin{example}\label{e:moebius}
(a) The Moebius contest \url{http://www.moebiuscontest.ru} is a prize for scientific papers.
Although this is not a peer review journal, this award plays an important role in setting reliability standards for young researches.
So Remark \ref{r:qual} is applicable.
Papers that are not reliable references do receive this prize, see (b) and \cite{Skm}.
So it is important either to make disclaimer to this effect, or to remove the word `scientific'
from the official information on the contest (\url{http://www.moebiuscontest.ru}).

%The reliability standards of the jury of the Moebius contest\linebreak
%(\url{http://www.moebiuscontest.ru}) are illustrated by the following.

(b) In  \cite[\S7]{Skm} I present the referee reports to \cite{Ak07} and to the earlier versions of \cite{Av14, Ru14, Go16} cited in \cite{Skm} (two of the four reports are in English).
I suggested to award these versions, making disclaimers that they are not reliable references.
However, the jury awarded them without such disclaimers.
For the version of \cite{Go16} the jury even publicly stated that it found the proof in that version complete
(as a jury member I voted against this statement).
A. Akopian was very decent to inform the jury several months after the award that he found a mistake in \cite{Ak07}; no update of \cite{Ak07} is publicly available.
The referees' remarks were later taken into account in \cite{Av14, Ru14} but not in \cite{Go16}.

%(c)

\end{example}

%%%%%%%%%%%%%%%%%%%%%%%%%%

\section{Appendix: letters to the Editors}\label{s:letters}

In this section we present our letters to the Editors concerning referee reports whose incompetence is justified in \S\ref{s:isrask17}.
First, this shows that before calling community's attention to this incompetence we spent quite some efforts to resolve the problem directly with the Editors (so that their willing to act upon incompetent reports is consistent policy, not momentary lapse).
Second, the letters (although written not by native English speakers) could be helpful for other authors to draft analogous letters in analogous situations.

\small

Recall the conventions from the beginning of \S\ref{s:peerex}.

\begin{remark}[to the Editors of Israel J. Math. on \cite{Sk17-2}]\label{r:ijm1}
{\bf (a; 29.01.2021)}

Dear Editors,

Hope you are fine and healthy.

Below [added in 2022: in (A)-(G) of Example \ref{r:revrevi}] please find my reply to the (attached) referee report on the manuscript

`Eliminating higher-multiplicity intersections in the metastable dimension range'

rejected from Israel Journal of Mathematics (IJM).

This reply shows that the report is incompetent.
Remark 5.3 of \cite{Sk17-3} shows that the report misuses the anonymous peer review system to promote an opinion which failed to be substantiated in an open discussion.
So could you please consider if the following would allow to maintain high reputation of IJM (see arXiv:2101.03745v1, Remark 2.2.Sk17 and \url{https://scirev.org/}):

(1) abolishing the rejection decision;

(2) sending the manuscript (whether the version submitted in 2020 or the attached one [added in 2022: \cite{Sk17-3}])
%the version attached to this letter
to an alternative referee;

(3a) sending my reply to the referee and asking him/her to present a revised report taking into account critical remarks presented in my reply
(the new version of the report can be written either on the version submitted in 2020 or on the attached one);

(3b) disregarding the referee report in case he/she would not provide such a revised report within time range corresponding to the standards of IJM;

(3c) if the Editors consider misusing the anonymous peer review system to promote an opinion which does not
stand an open discussion to be offence serious enough to disregard the referee report,
then I do not insist on (3a) and (3b).

%So I would be glad if the Editors of IJM send  to the same referee the revision and responses.
%A discussion of my reply by the Editors of IJM and sending them to the referee could maintain
%high standards of the peer review. However,

I do understand that moderation of a potential dispute between the referee and the author could be outside of purposes of the Editors.
So observe that most part of this dispute is irrelevant to publication of the submitted manuscript, see Example \ref{r:revrevi}.B below [added later: in this version, above].

I am open to compromise, e.g. to treat the attachment as a new submission if this would not result in delay of its refereeing.

Perhaps I should inform you that I am completing analogous reply to the referee report on the paper
\cite{AKS-1} rejected from IJM.
In that reply I justify that the report misuses the anonymous peer review system to promote an opinion which failed to be substantiated in an open discussion.

Sincerely yours, Arkadiy Skopenkov, \url{https://users.mccme.ru/skopenko/}.

\smallskip
{\bf (b; 26.03.2021)}  Dear Editors,

Hope you are fine and healthy.

Thank you for your letter.
Could you send me the referee report upon which your rejection decision of March 18 is based?

Recall that your previous rejection decision is based on an incompetent report that misuses the anonymous peer review system to promote an opinion which failed to be substantiated in an open discussion.
This judgement is justified in my letter of January 29 and is not questioned in your reply of March 18.
In your letter of March 18 you did not state that your previous rejection decision is retracted
(so that your new rejection decision is based only on the new report).
Could you please consider if making such a public statement would allow maintaining a high reputation of IJM
(see arXiv:2101.03745v1, Remark 2.2.Sk17 and \url{https://scirev.org/}).
For the math community it is important to know whether the anonymous peer review system is misused in IJM (cf. arXiv:2101.03745v1, Remarks 2.2.Sk16 and 2.2.Sk18 on Russian Math. Surveys [added in 2022: and on] Arnold J. Math.).
This is important independently of the journal's final decision on a particular paper.

Sincerely yours, Arkadiy Skopenkov, \url{https://users.mccme.ru/skopenko/}.

\smallskip
{\bf (c) 17.03.2022, to Michael Temkin, Editor in chief, and the Editors}

Shalom  Michael, Shalom Editors,

Hope you are fine and healthy.

Michael's letter of 29.03.2021 is delightful in a personal sense, because it sincerely attempted to resolve the problem in a friendly way, and Michael spent his valuable time to convince more experts to give their opinion.

Michael's letter of 29.03.2021 is disappointing in a professional sense, because it shows that Israel J. Math violates important principles of scientific discussions exposed
%does not maintain peer review standards.
% (as I suppose they are commonly understood).
in Remarks
%4.1, 4.2, 4.3 of  arXiv:2101.03745.
\ref{r:qual}, \ref{r:prin}, \ref{r:sugg}.
This is  justified in (a,b) and in [added later: (A)-(G) of] Example \ref{r:revrevi}.
%\S6 of  arXiv:2101.03745.
In particular, the more experts confirm the editorial decision based on an incompetent referee report without producing a new report sent to the author, the farther is the journal from a peer reviewed journal.

Your reaction suggested in Example \ref{r:viol}.a
%6.1.a
would be important.
E.g. you can

$\bullet$ restore the status of a peer reviewed journal (as suggested in Remark \ref{r:sugg}),
%4.3 of arXiv:2101.03745),
or

$\bullet$ publicly explain why the Editors disagree with specific parts of Remarks
%4.1, 4.2, 4.3 of  arXiv:2101.03745,
\ref{r:qual}, \ref{r:prin}, \ref{r:sugg},
or

$\bullet$ announce at the web page of the journal {\it `that the journal is not a peer reviewed journal, and so should not be used as a reliable reference and/or for jobs or grants evaluation'} (see Remark \ref{r:qual};
%4.1 of  arXiv:2101.03745;
this can be done implicitly by presenting no public reply to this letter).

Please find the pdf of this message at \url{https://www.mccme.ru/circles/oim/rese_inte.pdf}.
This message is a part of transparent public discussion as described in Remark \ref{r:transp}.c.
%4.4.c.
So in order to avoid confusion, could you read that remark before replying to this letter.

Sincerely Yours, Arkadiy Skopenkov.
\end{remark}

%Hopefully Israel J. Math. would choose the first option.
%I.e., to take decisions based on open discussion, not on justifications unavailable to the author
%(and to the math community), and so potentially incompetent / biased.
%coming from biased experts (as explained in \cite[Remark 5.3]{Sk16}).

\begin{remark}[to the Editors of Israel J. Math. on \cite{AKS-1}]\label{r:ijm-ks}
{\bf (a) 17.03.2022}

Dear Editors,

Hope you are fine and healthy.

Please find in Example \ref{r:aks-ijm}
%5.7 of arXiv:2101.03745
justification that the paper \cite{AKS-1}
%(cited there)
was rejected from  Israel J. Math. on the basis of an incompetent referee report.
In our opinion, this violates important principles of scientific discussions recalled in Remark \ref{r:prin}.ab.
%4.2.ab of arXiv:2101.03745.
Your reaction to this justification would be important, e.g. you can:

$\bullet$ restore the status of a peer reviewed journal (as suggested in Remark
%4.3.b of arXiv:2101.03745),
\ref{r:sugg}),
or

$\bullet$ publicly explain why the Editors disagree with specific parts of Remarks
%4.1, 4.2, and 4.3 of arXiv:2101.03745,
\ref{r:qual}, \ref{r:prin}, \ref{r:sugg},
or

$\bullet$ announce at the web page of the journal {\it `that the journal is not a peer reviewed journal, and so should not be used as a reliable reference and/or for jobs or grants evaluation'} (see Remark
%4.1 of arXiv:2101.03745;
\ref{r:qual}; this can be done implicitly by presenting no public reply to this letter).

Cf. Examples \ref{r:viol}.ab.
%6.1.abc of arXiv:2101.03745.

Please find the pdf of this message at \url{https://www.mccme.ru/circles/oim/rese_inte.pdf}.

Sincerely Yours, Roman Karasev and Arkadiy Skopenkov.

\smallskip
{\bf (b) (A reply of 18.03.2022 to a private letter from M. Temkin of 17.03.2022)}

Dear Michael,

Could you provide your public replies to  Arkadiy's, and/or to Roman and Arkadiy's messages of March 17,
so that your replies can be presented in an update of arXiv:2101.03745?
See Remark \ref{r:transp}.c.
%4.4.c of arXiv:2101.03745.

If not then  \url{https://www.mccme.ru/circles/oim/rese_inte.pdf} (and eventually arXiv:2101.03745)
will state (as a part of the public discussion) that no public reply is available.
This will be as informative to the math community as a public reply.

Sincerely Yours, Roman Karasev and Arkadiy Skopenkov.
\end{remark}

\begin{remark}\label{r:akspamsl}
(to the Editors I. Novik, D. Futer, J. Wang,  A. Folsom of Proc. of the Amer. Math. Soc. on \cite{AKS-2})

%Cc David Futer, dfuter@temple.edu Jiaping Wang jiaping@umn.edu,  Amanda Folsom afolsom@amherst.edu
%Dear Amanda, Dear Jiaping, Dear David,

{\bf (a) 20.03.2022}

Dear Isabella,

Hope you are fine and healthy.

Please find in Example \ref{r:akspamsd} our response to report X on our paper
\cite{AKS-2}
%arXiv:1908.08731v2
recently rejected from Proceedings of the AMS.
{\it There we justify that the report is incompetent, and violates important principles of scientific discussion recalled in Remark \ref{r:prin}.ab.}
%4.2.ab of arXiv:2101.03745.}

Let us quote the `Conclusion and recommendation' part of report Z:
{\it `I enjoyed reading this paper. It addresses important questions in a currently
very active area of topological and geometric combinatorics. The results are
non-trivial and interesting and the technique relatively new.
I gladly recommend this paper for publication in the Proceedings of the
American Mathematical Society.'}

Let us quote the last concluding paragraph of quick opinion Y:
{\it `I recommend proper refereeing. Given that there is a wealth of literature for degrees of equivariant maps,
a more careful verification is necessary to determine whether the main ingredient constitutes a new result. If this is the case, then this manuscript is of significant interest.'}

{\it Hence retraction of the rejection decision, and reconsideration of the paper, in our opinion, would allow to keep high reputation of Proceedings of the AMS.}
This can be done along the lines of Remark \ref{r:sugg}.b.
%4.3.b of arXiv:2101.03745.
{\it Could you please let us know if you plan to do that?}

{\it Reviews of handling A. Skopenkov's papers are published in \S\ref{s:peerex}, \S\ref{s:isrask17}.}
%\S5, \S6 of arXiv:2101.03745}.
E.g. for handling the previous version of the same paper see Examples \ref{r:viol}.a and \ref{r:aks-ijm}.
%5.7 of arXiv:2101.03745.

Best Regards, Roman Karasev and Arkadiy Skopenkov.

\smallskip
{\bf (b) 26.03.2022}

Dear Editors,

We write this letter in great respect to the American Mathematical Society.
This respect involves our belief that no AMS journal is an instrument of redistribution of jobs and grants in a way obstructing the progress of science,
so that mistakes in handling submitted papers are corrected, see Remark \ref{r:sugg}.
%4.3 of arXiv:2101.03745

Our letter of 20.03.2022 contains (in (b)) justification that the paper \cite{AKS-2}
was rejected from Proc. of the Amer. Math. Soc. on the basis of an incompetent referee report.
We thank Isabella for her reply.
Isabella's reply does not question that justification.
In our opinion, rejection of a paper on the basis of an incompetent report
violates important principles of scientific discussions recalled in Remark \ref{r:prin}.ab.
%4.2.ab of arXiv:2101.03745.
More importantly, leaving such a flaw uncorrected would show that a journal is not a peer review journal.
Isabella's reply does not demonstrate understanding of this problem,
and does not suggest or accept any way of resolving the problem.
(We did not object to rejection based on one {\it negative} report out of three, we objected to rejection based on an {\it incompetent} report.)
This calls for a transparent public discussion as described in Remark \ref{r:transp}.c.
%4.4.c of arXiv:2101.03745.
So your {\it public} reply to our letters would allow to keep high reputation of  Proc. of the Amer. Math. Soc.
In order to avoid confusion, could you read Remark \ref{r:transp}.c
%4.4.c of arXiv:2101.03745
before replying.

You can, in particular:

$\bullet$ carry public discussion of our justification (b) that report X is incompetent, as suggested in Remark \ref{r:sugg}
%4.3.b of arXiv:2101.03745
(and thus
%restore
confirm the status of a peer reviewed journal),
or

$\bullet$ publicly explain why the Editors disagree with specific parts of Remarks
%4.1, 4.2, and 4.3 of arXiv:2101.03745,
\ref{r:qual}, \ref{r:prin}, \ref{r:sugg},
or

$\bullet$ announce at the web page of the journal {\it `that the journal is not a peer reviewed journal, and so should
not be used as a reliable reference and/or for jobs or grants evaluation'} (see Remark
%4.1 of arXiv:2101.03745;
\ref{r:qual}; this can be done implicitly by presenting no public reply to this letter).

Cf. Examples \ref{r:skrms}, \ref{r:csmmj}, \ref{r:skajm}, \ref{r:ps-dcg}, and \ref{r:viol}.
%5.2, 5.4, 5.5, 5.6, 6.1 of arXiv:2101.03745.

Sincerely Yours, Roman Karasev and Arkadiy Skopenkov.
\end{remark}

%Emphasize???

%Unfortunately Isabella ignored our objection to taking the editorial decision upon an {\it incompetent} report)
%[only in a private letter: Thus Isabella's reply is the `straw man' logical fallacy.]
%This raises the question of whether Proc. of the Amer. Math. Soc. is a peer reviewed journal.

%No she agreed to the suggestion of Remark \ref{r:sugg}.b involving sending our justification to the referee for his/her potential objections.
%So we have to work under assumption that the report is indeed incompetent.
%Isabella explained that the journal rejects submissions that received at least one {\it negative} report, although we did not object to  that rule.
%Unfortunately Isabella ignored our objection to taking the editorial decision upon an {\it incompetent} report.

\begin{remark}[to the Editor H. Toru\'nczyk of Fund. Math. on \cite{KS21e}]\label{r:ksfm-l}

\textbf{(a) 27.12.2023}

Dear Henryk,

Thank you for your letters.
Could you send me the referee report(s) on which the Editorial decision is based?
If yes, could you allow me to publish them analogously to \cite{Rep}?

%\url{http://www.map.mpim-bonn.mpg.de/Template:Karasev\%26Skopenkov2020} ?

Hopefully the referee reports would

(1) explicitly name NEW results of the paper (the result stated in the abstract and restated as Corollary 1.3.3a,
Corollary 1.2.1b, Theorem 1.2.3, Corollaries 1.3.2ab, Corollary 1.3.5, and
the case k=2 of 1.2.1a, 1.3.3b, 1.2.2, 1.3.4, 1.3.7c);

(2) explicitly state that the results named in (1) do not meet the very high standards of Fundamenta Mathematicae;

(3) explicitly state where results of the removed Appendix of \cite{KS21e} are used in the proofs of the submitted paper.

This done, I will actually like the reports and the editorial decision made upon such clever reports.

[Here the text of Example \ref{r:ksfm}.a was presented]

Best Regards, Arkadiy.

PS If the reports and editorial decision are based on a confusion, then
we are ready to revise the paper emphasizing the distinction between old and new results.
In some cases this would make the paper less readable, so we need a referee's opinion
on whether he / she wants us to repeat several times that
for $k>2$ a result is known, and for $k=2$ (2-complexes in 4-manifolds) a result is new
(instead of having all this in compact form in Remark 1.3.7.a).

\smallskip
\textbf{(b) 15.01.2024}

Dear Henryk,

Hope you are fine and healthy.

Sending a report takes less than a minute.
So one has to conclude that there are no reports on the paper
which are up to the high standards of Fundamenta Mathematicae.
(Could you send me the reports without delay if the conclusion is wrong?)
So one has to conclude that the status of the paper changed back to `under review' (instead of `rejected').
Could you please confirm that (or send me the reports without delay if the conclusion is wrong) ?

Best regards, Arkadiy.

\smallskip
\textbf{\bf (c) 12.02.2024}

Dear Editors, Dear Henryk,

Hope you are fine and healthy.

I am sorry I did not receive any reply to the Dec 27, 2023 letter and the Jan 15, 2024 reminder.
Unfortunately, this leaves for me only to add to \S\ref{s:isrask17} of arXiv:2101.03745 [added later: of this paper] the following information:

[Here the text of Example \ref{r:ksfm}, before (a), was presented]

%[presented in (a); only the sentence involving reference to (b,c,d,e) is slightly updated].

However, I value collaboration with Polish colleagues very much.
So I decided to patiently send this reminder, and to wait another week waiting for the proper refereeing process.
That is, for either sending me the referee reports, or
confirming that the paper (submitted in April 2022) is back at the `under review' status,
and setting a timeframe for presenting the (new) reports together with Editorial decision based on these reports.

Best Regards, Arkadiy.

%[Repetition of the above information in this further polite reminder is omitted]

\smallskip
\textbf{(d) 27.02.2024}

Dear Editors, Dear Henryk,
%Dear J'ozef, Dear James,  

However reluctantly, I updated  \cite{Sk21d'} adding Example \ref{r:ksfm}.
It justifies that the paper arXiv:2112.06636 (shortened by deleting sections 2.4.and 2.5)
was rejected from Fundamenta Mathematicae because of a logical fallacy in the Editors' letter, and a wrong statement in the Editors' letter.
(I am sorry I did not receive any reply to the Dec 27, 2023 letter and the Jan 15, 2024 reminder.
I am sorry I did not receive the information I asked for in my Feb 12, 2024 letter.)

I value collaboration with Polish colleagues very much.
So I postponed updating Example \ref{r:viol} of \cite{Sk21d'}, and updating arXiv:2101.03745.
I am patiently sending this reminder, and waiting another week for confirmation of the proper refereeing process. 
That is, for either sending me the referee reports, or confirming that the paper (submitted in April 2022) is back at the `under review' status, and setting a timeframe for presenting the (new) reports together with Editorial decision based on these reports.
I would be glad to either add positive updated information to Example \ref{r:ksfm}, or to leave (upon your wish) only positive updated information in Example \ref{r:ksfm} (in both cases moving Example \ref{r:ksfm} to a different section). 

Please observe that I am asking for the above confirmation of the proper refereeing process, which could be done immediately.
I am not asking for immediate sending me updated referee reports, which perhaps requires more than 2 months passed since my Dec 27, 2023 letter (although the paper is submitted in April, 2022).  

Sincerely Yours, Arkadiy Skopenkov.

PS Mathematicians (both Polish and foreign ones) would be grateful if you allow me to publish the Editors' Dec 8, 2023 letter giving grounds for rejection (perhaps except the first and the last paragraphs not related to such grounds).
This will help them to make their own judgement whether my criticism of that letter is just.
Lack of your agreement to publish that letter would be a serious confirmation that my criticism is just.
\end{remark}

\begin{remark}[to the Principal Editor R. Kenyon of Intern. Math. Res. Notes on \cite{DS-3, DS-5}] \label{r:dsirmn-l}
\textbf{(a) May 6, 2024}

Dear Principal Editor Rick Kenyon, Dear Editors

Hope you are fine and healthy.

Thank you for handling our paper ID IMRN-2024-227, A QUADRATIC ESTIMATION FOR THE K\"UHNEL CONJECTURE ON EMBEDDINGS.
The submitted version is obtained from \cite{DS-3} by deleting \S3, \S4 and footnotes.
In Example \ref{r:dsirmn}.abc we justify that the referee report on the paper is incompetent.
Most of the justification is accessible to mathematicians not specialized in the area.
Hence in our opinion, retraction of the rejection decision, and reconsideration of the paper,
would allow to keep the high reputation of the journal, cf. Remarks \ref{r:qual} and \ref{r:sugg}.a.
%4.1 and 4.3.a of arXiv:2101.03745.
This can be done by sending our response to the referee and / or assigning a different referee
(see detailed description of what we consider to be the procedure of quality research journals in Remark \ref{r:sugg}.b).
%of arXiv:2101.03745).
Could you please let us know if you plan to do that?

This discussion can be carried on the initially submitted version, or on the revision \cite{DS-5}
(although the referee report is incompetent, we used it to slightly improve the introduction).

We would also be grateful if you could check if the referee has no conflict of interest in refereeing our paper.

Reviews of handling A. Skopenkov’s papers are published at arXiv:2101.03745 [added later: in this paper], \S\ref{s:peerex}, \S\ref{s:isrask17}.
We would be glad to add there a positive review of handling this paper in IMRN, no matter if the paper is accepted or rejected.

Best regards, Arkadiy Skopenkov and Slava Dzhenzher.

[Here the text of Example \ref{r:dsirmn}.abc was presented]

\smallskip
\textbf{(b) May 20, 2024}

Dear Principal Editor Rick Kenyon, Dear Editors, 

Hope you are fine and healthy.

We are sorry we did not receive any reply to our May 6, 2024 letter.
Unfortunately, this leaves for us only to add to \S\ref{s:isrask17} of
arXiv:2101.03745 [added later: of this paper] Example \ref{r:dsirmn},
%6.10 from the attachment,
and to make corresponding additions on IMRN in Example \ref{r:viol}.a there.  

However, we value collaboration with IMRN very much.
So we decided to patiently send this reminder, and to wait another week waiting for your agreement to discuss our response to the referee report with a referee (old or new).
We are asking if you could confirm your agreement soon; we do understand that preparation of a new report addressing our response may require some time. 

Please let us know if technically submission of the revision is more convenient for you, to handle the situation involving an incompetent referee report.  

Best Regards, Arkadiy Skopenkov and Slava Dzhenzher.

\smallskip
\textbf{(c) April 7, 2025} (cf. footnote \ref{f:dsirmn})

Dear Principal Editor Rick Kenyon, Dear Editors,

Hope you are fine and healthy.

Thank you for handling our paper [DS5]:=ID IMRN-2024-608,
A QUADRATIC ESTIMATION FOR THE K\"UHNEL CONJECTURE ON EMBEDDINGS.
The previous version [DS3]:=ID IMRN-2024-227 submitted in February 2024
was rejected in March 2024 upon an incompetent referee report
(see our justification of incompetence in
%our 6 May 2024 letter
Example \ref{r:dsirmn}.abc).
The version [DS5] resubmitted in May 2024 is obtained from arXiv:2208.04188v5 by deleting \S3, \S4 and footnotes.

[Here the text of Example \ref{r:dsirmn5}.a was presented]

Hence, in our opinion, retraction of the rejection decision and reconsideration of the paper
would allow keeping the high reputation of the journal, cf. Remarks \ref{r:qual} and \ref{r:sugg}.a.
%4.1 and 4.3.a of arXiv:2101.03745.
This can be done by sending our responses to the referees and/or assigning different referees
(see detailed description of what we consider to be the procedure of quality research journals in Remark \ref{r:sugg}.b).
%4.3.b of arXiv:2101.03745).
Could you please let us know if you plan to do that?

This discussion can be carried upon the resubmitted version [DS5].
Alternatively, we can resubmit a shortened version of arXiv:2208.04188v7
(although the report \#1 is incompetent, and the report \#2 cannot contribute
to the rejection decision, we used the reports to improve the presentation).

Reviews of handling A. Skopenkov’s papers are published at arXiv:2101.03745 [added later: in this paper],
%(and at \linebreak \url{https://old.mccme.ru//circles//oim/rese_inte.pdf}),
\S\ref{s:peerex}, \S\ref{s:isrask17}.
We would be glad to add there a positive review of handling this paper in IMRN,
no matter if the paper is accepted or rejected.

For your convenience, please find the pdf of this letter in the attachment.

Best regards, Slava Dzhenzher and Arkadiy Skopenkov.

[Here the text of Example \ref{r:dsirmn5}.(b)-(d6) was presented, adding to the end of (b) the phrase
`At the beginning of each item a quotation from referee report is given in italics', and omitting quotation marks at the beginning of each item.]

\smallskip
\textbf{(d) April 13, 2025}

Dear Secretary Mrs  Emma-Louise Smith, Dear Principal Editor Rick Kenyon, Dear Editors,

Many thanks for your quick same-day replies.
The replies contradict each other.
Unless instructed by you otherwise, we have to act upon the reply of Rick Kenyon.

The reply is very interesting to the math community, as a revelation of the review integrity standards of IMRN.
So would you kindly allow us to publish the reply?
Technically speaking, would you release the copyright for the reply to the public domain?
Please feel free to make any changes in the reply before making it public.

Sincerely Yours, Arkadiy Skopenkov and Slava Dzhenzher.

PS Unfortunately, the reply of Rick Kenyon leaves for us only to add to \cite{Sk21d'} (and soon to arXiv:2101.03745) our critical remarks on the referee report, our letters to the Editors, and some conclusions on IMRN (see Examples \ref{r:dsirmn}-\ref{r:dsirmn5}, Remark \ref{r:dsirmn-l}, and Example \ref{r:viol}.a,
%Examples 6.7-6.8, Remark 7.5, and Example 6.1.a,
respectively).
Those conclusions are justified by the reply of Rick Kenyon.
We are glad that they are not justified by the reply of Emma-Louise Smith
(unfortunately, we cannot act upon that reply).

However, we value collaboration with IMRN very much.
So we decided to patiently send this request for a public reply, and to wait another several days before updating  arXiv:2101.03745.

If your new reply would keep high reputation of the journal by retraction of the rejection decision,
and reconsideration of our paper, then we remove from \cite{Sk21d'}
%\url{https://old.mccme.ru//circles//oim/rese_inte.pdf}
(and from an update of arXiv:2101.03745)
Examples \ref{r:dsirmn}-\ref{r:dsirmn5}, Remark \ref{r:dsirmn-l}, and the information on IMRN from Example \ref{r:viol}.a.
Later they will be replaced by the final information (hopefully praising the review integrity standards of IMRN, no matter if our paper is accepted or rejected).
Please let us know if technically submission of the revision is more convenient for you,
to handle the situation involving an incompetent referee report.

Otherwise in an update of  arXiv:2101.03745 we shall either add your reply,
or the information that you did not allow us to make it public.
The first option would be a sound confirmation that the Editors are willing to bear responsibility for those grounds that motivated their decision to ignore our (openly published) critical remarks on the referee report.
We are afraid that the second option would be a sound confirmation that our (openly published) critical remarks on the referee report are justified, so that  openly ignoring them will jeopardize the reputation of the journal.
See arXiv:2101.03745v3, Remark \ref{r:transp}.c on transparency.
%4.4.c.
\end{remark}

{\it In this list by stars I marked books, surveys, expository papers, as well as internet materials which I consider accessible to non-specialists.}

\end{document}